\documentclass[aos,seceqn,citesort,dvips]{arximspdf}
\usepackage{multirow}
\usepackage{graphicx}

\doi{10.1214/09-AOS780}
\volume{38}
\issue{4}
\pubyear{2010}
\firstpage{2092}
\lastpage{2117}

\makeatletter

\newcommand{\rmP}{\mathrm{P}}
\newcommand{\iint}{\int\!\!\int}

\newtheorem{theorem}{Theorem}[section]
\newtheorem{lemma}{Lemma}[section]
\newproclaim{Step}{Step}

\makeatother

\begin{document}
\begin{frontmatter}

\title{Penalized variable selection procedure for Cox models with
semiparametric relative risk}
\runtitle{Variable selection for semiparametric Cox models}

\begin{aug}
\author[A]{\fnms{Pang} \snm{Du}\corref{}\ead[label=e1]{pangdu@vt.edu}},
\author[B]{\fnms{Shuangge} \snm{Ma}\ead[label=e2]{shuangge.ma@yale.edu}} and
\author[C]{\fnms{Hua} \snm{Liang}\thanksref{t2}\ead[label=e3]{hliang@bst.rochester.edu}}
\runauthor{P. Du, S. Ma and H. Liang}
\affiliation{Virginia Tech, Yale University and University of Rochester}
\address[A]{P. Du\\
Department of Statistics\\
Virginia Tech\\
Blacksburg, Virginia 24061\\
USA\\
\printead{e1}}
\address[B]{S. Ma\\
Department of Epidemiology\\
\quad and Public Health\\
Yale University\\
School of Medicine\\
New Haven, Connecticut 06520\\
USA\\
\printead{e2}}
\address[C]{H. Liang\\
Department of Biostatistics\\
\quad and Computational Biology\\
University of Rochester\\
Rochester, New York 14642\\
USA\\
\printead{e3}}
\end{aug}

\thankstext{t2}{Supported in part by
NIH/NIAID Grant AI59773 and NSF Grant DMS-08-06097.}

\received{\smonth{8} \syear{2009}}
\revised{\smonth{12} \syear{2009}}

%
\begin{abstract}
We study the Cox models with semiparametric relative risk, which can be
partially linear with one nonparametric component, or multiple additive
or nonadditive nonparametric components. A penalized partial likelihood
procedure is proposed to simultaneously estimate the parameters and
select variables for both the parametric and the nonparametric parts.
Two penalties are applied sequentially. The first penalty, governing
the smoothness of the multivariate nonlinear covariate effect function,
provides a smoothing spline ANOVA framework that is exploited to derive
an empirical model selection tool for the nonparametric part. The
second penalty, either the smoothly-clipped-absolute-deviation (SCAD)
penalty or the adaptive LASSO penalty, achieves variable selection in
the parametric part. We show that the resulting estimator of the
parametric part possesses the oracle property, and that the estimator
of the nonparametric part achieves the optimal rate of convergence. The
proposed procedures are shown to work well in simulation experiments,
and then applied to a real data example on sexually transmitted
diseases.
\end{abstract}

%
\begin{keyword}[class=AMS]
\kwd[Primary ]{62N01}
\kwd{62N03}
\kwd[; secondary ]{62N02}.
\end{keyword}
\begin{keyword}
\kwd{Backfitting}
\kwd{partially linear models}
\kwd{penalized variable selection}
\kwd{proportional hazards}
\kwd{penalized partial likelihood}
\kwd{smoothing spline ANOVA}.
\end{keyword}

\end{frontmatter}

\section{Introduction}\label{sec:intr}

In survival analysis, a problem of interest is to identify relevant
risk factors and evaluate their contributions to survival time. Cox
proportional hazards (PH) model is a popular approach to study the
influence of covariates on survival outcome. Conventional PH models
assume that covariates have a log-linear effect on the hazard function.
These PH models have been studied by numerous authors; see, for
example, the references in \cite{kal02}. The log-linear assumption
can be too rigid in practice, especially when continuous covariates are
present. This limitation motivates PH models with nonparametric
relative risk. Some examples are
\cite{zucker90,osu93,fan97,huang00,huang06}. However,
nonparametric models may suffer from the curse of dimensionality. They
also lack the easy interpretation in parametric risk models. PH models
with semiparametric relative risk strike a good balance by allowing
nonparametric risk for some covariates and parametric risk for others.
The benefits of such models are two-folds. First, they have the merits
of models with parametric risk, including easy interpretation, easy
estimation and easy inference. Second, their nonparametric part allows
a flexible form for some continuous covariates whose patterns are
unexplored and whose contribution cannot be assessed by simple
parametric models. For example, \cite{huang99} proposed efficient
estimation for a partially linear Cox model with additive nonlinear
covariate effects. Reference \cite{cai07} studied partially linear hazard
regression for multivariate survival data with time-dependent
covariates via a profile pseudo-partial likelihood approach, where the
only nonlinear covariate effect was estimated by local polynomials. But
these models are limited to one nonparametric component or additive
nonparametric components, ignoring the possible interactions between
different nonparametric components. Reference \cite{yin08} proposed a partially
linear additive hazard model whose nonlinear varying coefficients
represent the interaction between the time-dependent nonlinear
covariate and other covariates.

Variable selection in survival data has drawn much attention in the
past decade. Traditional procedures such as Akaike information
criterion (AIC) and Bayesian information criterion (BIC), as noted by
\cite{breiman96}, suffer from the lack of stability and lack of
incorporating stochastic errors inherited in the stage of variable
selection. References \cite{tibs97} and \cite{zou08} extended,
respectively, the LASSO and the adaptive LASSO variable selection
procedures to the Cox model. Reference \cite{fanli02} extended the
nonconcave penalized likelihood approach \cite{scad01} to the Cox PH
models. Reference \cite{cai05} studied variable selection for multivariate
survival data. The Cox models considered in these three papers all
assumed a linear form of covariate effects in the relative risk. More
recently, \cite{johnson08} and \cite{jlz08} proposed procedures for
selecting variables in semiparametric linear regression models for
censored data, where the dependence of response over covariates was
also assumed to be of linear form. Hence, the aforementioned variable
selection procedures are limited in their rigid assumption of
parametric covariate effects which may not be realistic in practice. We
will fill in these gaps in three aspects: (i) our models are flexible
with semiparametric relative risk, which allows nonadditive
nonparametric components, without limiting to single or additive
nonlinear covariate effects; and (ii) our approach can simultaneously
estimate the parametric coefficient vector and select contributing
parametric components; and (iii) our approach also provides a model
selection tool for the nonparametric components.

Let the hazard function for a subject be
%
\begin{equation}\label{eq:mod}
h(t)=h_0(t)\exp[\bolds{\beta}^T U+\eta(W)],
\end{equation}
where $h_0$ is the unknown baseline hazard, $Z^{T}=(U^{T},W^{T})$ is
the covariate vector, $\bolds{\beta}$ is the unknown
coefficient vector, and $\eta(w)=\eta(w_1,\ldots,w_q)$ is an unknown
multivariate smooth function. We propose a doubly penalized profile
partial likelihood approach for estimation, following the general
profile likelihood framework set up by \cite{murphy00}. Given
$\bolds{\beta}$, $\eta$ is estimated by smoothing splines
through the minimization of a penalized log partial likelihood. Then
the smoothing spline ANOVA decomposition not only allows the natural
inclusion of interaction effects but also provides the basis for
deriving an empirical model selection tool. After substituting the
estimate of $\eta$, we obtain a profile partial likelihood, which is
then penalized to get an estimate of $\bolds{\beta}$. To
achieve variable selection in $\bolds{\beta}$, we use the
smoothly clipped absolute deviation (SCAD) penalty.
We show that our estimate of $\eta$ achieves the
optimal convergence rate, and our estimate of $\bolds{\beta}$
possesses the oracle property such that the true zero coefficients
are automatically estimated as zeros and the remaining coefficients
are estimated as well as if the correct submodel were known in
advance. Our numerical studies reveal that the proposed method is
promising in both estimation and variable selection. We
then apply it to a study on sexually transmitted diseases with 877
subjects.

The rest of the article is organized as follows.
Section \ref{sec:meth} gives the details of the proposed method, in
the order of model description and estimation procedure
(Section \ref{ssec:est}), model selection in the nonparametric part
(Section \ref{ssec:npc}), asymptotic properties (Section \ref
{ssec:asym}), and
miscellaneous issues (Section \ref{ssec:misc}) like standard error
estimates and smoothing parameter selection. Section \ref{sec:num}
presents the empirical studies, and Section \ref{sec:app} gives an
application study. Remarks in Section \ref{sec:disc} conclude the
article.

\section{Method}\label{sec:meth}

Let $T$ be the failure time and $C$ be the right-censoring time.
Assume that $T$ and $C$ are conditionally independent given the
covariate. The observable random variable is
$(X,\Delta,Z)$, where $X=\min(T,C)$, $\Delta=I_{[T\le C]}$, and
$Z=(U,W)$ is the covariate vector with $U\in\mathbb{R}^{d}$ and
$W\in\mathbb{R}^{q}$. With $n$ i.i.d. $(X_{i},\Delta_{i},Z_{i}),
i=1,\ldots,n$, we assume a Cox model for the hazard function as
in~(\ref{eq:mod}).

\subsection{Estimation and variable selection for parametric parts}
\label{ssec:est}

Let $Y_{i}(t)=I_{[X_{i}\ge{t}]}$. We propose to estimate
$(\bolds{\beta},\eta)$ through a penalized profile partial
likelihood approach. Given $\bolds{\beta}$, $\eta$ is
estimated as the minimizer of the penalized partial likelihood
%
\begin{eqnarray}\label{equ:pplk1}\qquad\quad
l_{\bolds{\beta}}(\eta)&\equiv&-\frac{1}{n}\sum_{i=1}^{n}
\Delta_{i} \Biggl\{
U_{i}^{T}\bolds{\beta}+\eta(W_{i})
-\log\sum_{k=1}^{n}Y_{k}(X_{i})\exp[U_{k}^{T}\bolds{\beta}
+\eta(W_{k})] \Biggr\}\nonumber\\[-9pt]\\[-9pt]
&&{}+\lambda{J}(\eta),\nonumber
\end{eqnarray}
where the summation is the negative log partial likelihood
representing the goodness-of-fit, $J(\eta)$ is a roughness penalty
specifying the smoothness of $\eta$, and $\lambda>0$ is a smoothing
parameter controlling the tradeoff. A popular choice for $J$ is the
$L_2$-penalty which yields tensor product cubic splines (see, e.g.,
\cite{gu02}) for multivariate $W$. Note that
$\eta$ in (\ref{equ:pplk1}) is identifiable up to a constant, so we
use the constraint $\int\eta=0$.

Once an estimate $\hat{\eta}$ of $\eta$ is obtained, the estimator
of $\bolds{\beta}$ is then the maximizer of the penalized
profile partial likelihood
%
\begin{eqnarray}\label{equ:pplk2}\qquad\quad
l_{\hat{\eta}}(\bolds{\beta})&\equiv&\sum_{i=1}^{n}\Delta_{i} \Biggl\{
U_{i}^{T}\bolds{\beta}+\hat{\eta}(W_{i})
-\log\sum_{k=1}^{n}Y_{k}(X_{i})\exp[U_{k}^{T}\bolds{\beta}
+\hat{\eta}(W_{k})] \Biggr\}\nonumber\\[-9pt]\\[-9pt]
&&{}-n\sum_{j=1}^{d}p_{\theta_{j}}(|\beta_{j}|),\nonumber
\end{eqnarray}
where $p_{\theta_{j}}(|\cdot|)$ is the SCAD penalty on $\bolds{\beta}$
\cite{scad01}.

The detailed algorithm for our estimation procedure is as
follows.

\begin{enumerate}[Step 2.]
\item[Step 1.]
Find a proper initial estimate $\hat{\beta}^{(0)}$. We note that, as
long as the initial estimate is reasonable, convergence to the true
optimizer can be achieved. Difference choices of the initial
estimate will affect the number of iterations needed but not the
convergence itself.

\item[Step 2.]
Let $\hat{\bolds{\beta}}{}^{(k-1)}$ be the estimate of
$\bolds{\beta}$ before the $k$th iteration. Plug
$\hat{\bolds{\beta}}{}^{(k-1)}$ into (\ref{equ:pplk1}) and solve
for $\eta$ by minimizing the penalized partial likelihood
$l_{\hat{\bolds{\beta}}{}^{(k-1)}}(\eta)$. Let $\hat{\eta}^{(k)}$
be the estimate thus obtained.

\item[Step 3.]
Plug $\hat{\eta}^{(k)}$ into (\ref{equ:pplk2}) and solve for
$\bolds{\beta}$ by maximizing the penalized profile partial
likelihood $l_{\hat{\eta}^{(k)}}(\bolds{\beta})$. Let
$\hat{\bolds{\beta}}{}^{(k)}$ be the estimate thus obtained.
\item[Step 4.]
Replace $\hat{\bolds{\beta}}{}^{(k-1)}$ in step 2 by $\hat{\bolds{\beta
}}{}^{(k)}$ and
repeat steps 2 and 3 until convergence to obtain the final estimates
$\hat{\bolds{\beta}}$ and $\hat{\eta}$.
\end{enumerate}
Our experience shows that the algorithm usually converges quickly
within a few iterations. As in the classical Cox proportional
hazards model, the estimation of baseline hazard function is of less
interest and not required in our estimation procedure.

In step 3, we use a one-step approximation to the SCAD penalty
\cite{zouli08}. It transforms the SCAD penalty problem to a
LASSO-type optimization, where the celebrated LARS algorithm
proposed in \cite{efron04} can be used. Let
$l_{\hat{\eta}}(\bolds{\beta})$ be the profile log partial
likelihood in step 3, and
$I(\bolds{\beta})=-\nabla^{2}l_{\hat{\eta}}(\bolds{\beta})$
be the Hessian matrix, where the derivative is with respect to
$\bolds{\beta}$ treating $\hat{\eta}$ as fixed. Compute the
Cholesky decomposition of $I(\hat{\bolds{\beta}}{}^{(k-1)})$ such
that $I(\hat{\bolds{\beta}}{}^{(k-1)})=V^{T}V$. Let
$A=\{j\dvtx p'_{\theta_{j}}(|\hat{\beta}_{j}^{(k-1)}|)=0\}$ and
$B=\{j\dvtx p'_{\theta_{j}}(|\hat{\beta}_{j}^{(k-1)}|)>0\}$. Decompose
$V$ and the new estimate $\hat{\bolds{\beta}}{}^{(k)}$
accordingly such that $V=[V_{A},V_{B}]$ and
$\hat{\bolds{\beta}}{}^{(k)}=
(\hat{\bolds{\beta}}{}^{(k)^{T}}_{A},
\hat{\bolds{\beta}}{}^{(k)^{T}}_{B})^{T}$.

\begin{enumerate}[(Step 3a)]
\item[(Step 3a)]
Let $y=V\hat{\bolds{\beta}}{}^{(k-1)}$. Then for each $j\in B$,
replace the $j$th column of $V$ by setting
$v_{j}=v_{j}\frac{\theta_{j}}{p_{\theta_{j}}'(|\hat{\beta}_{j}^{(k-1)}|)}$.
\item[(Step 3b)]
Let $H_{A}=V_{A}(V_{A}^{T}V_{A})^{-1}V_{A}^{T}$ be the projection
matrix to the column space of $V_{A}$. Compute
$y^{*}=y-H_{A}y$ and $V_{B}^{*}=V_{B}-H_{A}V_{B}$.
\item[(Step 3c)]
Apply the LARS algorithm to solve
\[
\hat{\bolds{\beta}}{}^{*}_{B}=
\mathop{\arg\min}_{\bolds{\beta}}\biggl\{\frac{1}{2}
\|y^{*}-V_{B}^{*}\bolds{\beta}\|^{2}+
n\sum_{j\in{B}}\theta_{j}|\beta_{j}|\biggr\}.
\]
\item[(Step 3d)]
Compute
$\hat{\bolds{\beta}}{}^{*}_{A}
=(V_{A}^{T}V_{A})^{-1}V_{A}^{T}(y-V_{B}\hat{\bolds{\beta}}{}^{*}_{B})$
to obtain
$\hat{\bolds{\beta}}{}^{*}=
(\hat{\bolds{\beta}}{}^{*{T}}_{A},
\hat{\bolds{\beta}}{}^{*{T}}_{B})^{T}$.
\item[(Step 3e)]
For $j\in{A}$, set
$\hat{\beta}_{j}^{(k)}=\hat{\beta}_{j}^{*}$.
For $j\in{B}$, set
$\hat{\beta}_{j}^{(k)}=\hat{\beta}_{j}^{*}
\frac{\theta_{j}}{p'_{\theta_{j}}(|\hat{\beta}_{j}^{(k-1)}|)}$.
\end{enumerate}

\subsection{Model selection for nonparametric component}
\label{ssec:npc}

While the SCAD\break penalty takes care of variable selection for the
parametric components, we still need an approach to assess the
structure of the nonparametric components. In this section, we will
first transform the profile partial likelihood problem in
(\ref{equ:pplk1}) to a density estimation problem with biased
sampling, and then derive a model selection tool based on the
Kullback--Leibler geometry. In this part, we treat
$\bolds{\beta}$ as fixed, taking the value from the previous
step in the algorithm.

Let $(i_{1},\ldots,i_{N})$ be the indices for the failed subjects.
Then the profile partial likelihood in (\ref{equ:pplk1}) for
estimating $\eta$ is
\[
\prod_{i=1}^{n} \biggl[\frac{e^{U_{i}^{T}\bolds{\beta}+\eta(W_{i})}}
{\sum_{k=1}^{n}Y_{k}(X_{i})e^{U_{k}^{T}\bolds{\beta}+
\eta(W_{k})}} \biggr]^{\Delta_{i}}
=\prod_{p=1}^{N} \biggl[\frac{e^{U_{i_p}^{T}\bolds{\beta}+
\eta(W_{i_p})}}
{\sum_{k=1}^{n}Y_{k}(X_{i_p})e^{U_{k}^{T}\bolds{\beta}+
\eta(W_{k})}} \biggr].
\]
Consider the empirical measure $\rmP_{n}^{w}$ on the discrete
domain $\mathcal{W}_{n}=\{W_{1},\ldots,\break W_{n}\}$ such that
$\int{f}\,d\rmP_{n}^{w}=\frac{1}{n}\sum_{i=1}^{n}f(W_{i})$. Then
$e^{\eta}/\int e^{\eta}\,d\rmP_{n}^{w}$ defines a density function
on $\mathcal{W}_{n}$. Let $a_{1}(\cdot), \ldots, a_{N}(\cdot)$ be
weight functions defined on the discrete domain $\mathcal{W}_{n}$
such that
$a_{p}(W_k)=Y_{k}(X_{i_p})e^{U_{k}^{T}\bolds{\beta}}$,
$p=1,\ldots,N$. Alternatively, one can think of $a_p$'s as vectors
of weights with length $n$. Then each term in the profile partial
likelihood, with the constant $n$ ignored, becomes
$a_{p}(W_{i_{p}})e^{\eta(W_{i_{p}})}/\int
a_{p}(w)\times e^{\eta(w)}\,d\rmP_{n}^{w}$. Thus, this resembles a
density estimation problem with bias introduced by the known weight
function $a_{p}(\cdot)$.

For two density estimates $\eta_{1}$ and $\eta_{2}$ in the above
pseudo biased sampling density estimation problem, define their
Kullback--Leibler distance as
%
\begin{eqnarray}\label{equ:ekl}
\operatorname{KL}(\eta_{1},\eta_{2})&=&\frac{1}{N}\sum_{p=1}^{N} \biggl\{
\frac{\int(\eta_{1}(w)-\eta_{2}(w))a_{p}(w)e^{\eta_{1}(w)}\,d\rmP_{n}^{w}}
{\int a_{p}(w)e^{\eta_{1}(w)}\,d\rmP_{n}^{w}}\nonumber\\
&&\hspace*{34.17pt}{}-\log\int{a}_{p}(w)e^{\eta_{1}(w)}\,d\rmP_{n}^{w}\\
&&\hspace*{34.17pt}\hspace*{44.69pt}{}+\log\int{a}_{p}(w)
e^{\eta_{2}(w)}\,d\rmP_{n}^{w} \biggr\}.\nonumber
\end{eqnarray}
Let $\eta_{0}$ be the true function. Suppose the estimation of
$\eta_{0}$ has been done in a space~$\mathcal{H}_{1}$, but in fact
$\eta_{0}\in\mathcal{H}_{2}\subset\mathcal{H}_{1}$. Let
$\hat{\eta}$ be the estimate of $\eta_{0}$ in $\mathcal{H}_{1}$. Let
$\tilde{\eta}$ be the Kullback--Leibler projection of $\hat{\eta}$ in
$\mathcal{H}_{2}$, that is, the minimizer of
$\operatorname{KL}(\hat{\eta},\eta)$ for $\eta\in\mathcal{H}_{2}$, and
$\eta_{c}$ be the estimate from the constant model. Set
$\eta=\tilde{\eta}+\alpha(\tilde{\eta}-\eta_{c})$ for $\alpha$ real.
Differentiating $\operatorname{KL}(\hat{\eta},\eta)$ with respect to
$\alpha$ and evaluating at $\alpha=0$, one has
\[
\sum_{p=1}^{N}\frac{\int(\tilde{\eta}(w)
-\eta_{c}(w))a_{p}(w)e^{\hat{\eta}(w)}\,d\rmP_{n}^{w}}
{\int a_{p}(w)e^{\hat{\eta}(w)}\,d\rmP_{n}^{w}}=
\sum_{p=1}^{N}\frac{\int(\tilde{\eta}(w)
-\eta_{c}(w))a_{p}(w)e^{\tilde{\eta}(w)}\,d\rmP_{n}^{w}}
{\int a_{p}(w)e^{\tilde{\eta}(w)}\,d\rmP_{n}^{w}},
\]
which, through straightforward calculation, yields
\[
\operatorname{KL}(\hat{\eta},\eta_{c})=
\operatorname{KL}(\hat{\eta},\tilde{\eta})+
\operatorname{KL}(\tilde{\eta},\eta_{c}).
\]
Hence, the ratio
$\operatorname{KL}(\hat{\eta},\tilde{\eta})/\operatorname{KL}(\hat{\eta
},\eta_{c})$
can be used to diagnose the feasibility of a reduced model
$\eta\in\mathcal{H}_{2}$: the smaller the ratio is, the more
feasible the reduced model is.

\subsection{Asymptotic results}
\label{ssec:asym}

Denote by $\mathcal{H}^{m}(\mathcal{W})$
the Sobolev space of functions on $\mathcal{W}$
whose $m$th order partial derivatives are square integrable.
Let
\[
\mathcal{H}= \biggl\{\eta\in\mathcal{H}^{m}(\mathcal{W}),
\int_{\mathcal{W}}\eta(w)\,dw=0 \biggr\},
\]
and $\hat{\eta}^*$
be the estimate of $\eta_0$ in $\mathcal{H}$
that minimizes the penalized partial likelihood
%
\begin{eqnarray}\label{eq:plk}\qquad
&&-\frac{1}{n}\sum_{i=1}^{n}\int_{\mathcal{T}} \Biggl\{
U_{i}^{T}\bolds{\beta}+\eta(W_{i})
-\log\sum_{k=1}^{n}Y_{k}(t)\exp[U_{k}^{T}\bolds{\beta}
+\eta(W_{k})] \Biggr\}\,dN_i(t)\nonumber\\[-8pt]\\[-8pt]
&&\qquad{}+\frac{\lambda}{2}{J}(\eta).\nonumber
\end{eqnarray}
Note that $\mathcal{H}$ is an infinite-dimensional function space.
Hence, in practice, the minimization of (\ref{eq:plk}) is usually
performed in a data-adaptive finite-dimensional space
\[
\mathcal{H}_{n}=\mathcal{N}_{J}\oplus\operatorname{span}
\{R_{J}(W_{i_l},\cdot)\dvtx l=1,\ldots,q_n \},
\]
where $\mathcal{N}_{J}=\{\eta\in\mathcal{H}, J(\eta)=0\}$ is the
null space of $J$, and $R_{J}$ is the \textit{reproducing
kernel} (see, e.g., \cite{wahba90}) in its complement space
$\mathcal{H}_{J}=\mathcal{H}\ominus\mathcal{N}_{J}$, and
$\{W_{i_1},\ldots,W_{i_{q_n}}\}$ is a random subset of $\{W_i\dvtx
i=1,\ldots,n\}$. When $q_n=n$, one selects all the
$W_i,i=1,\ldots,n$ as the knots. This is the number of knots used
in conventional smoothing splines. However, under the regression
setting, \cite{kg04} showed that a $q_n$ of the order
$n^{2/(r+1)+\varepsilon},\forall\varepsilon>0$ is sufficient to yield an
estimate with the optimal convergence rate. Here $r$ is a constant
associated with the Sobolev space $\mathcal{H}$, for example, $r=2m$ for
splines of order $m$ (one-dimension $w$) and
$r=2m-\delta,\forall\delta>0$ for tensor product splines
(multi-dimension $w$). We shall show that such an order for $q_n$
also works for the $\eta$ estimation in our partially linear Cox
model.

Let $s_n[f;\bolds{\beta},\eta](t)=\frac{1}{n}\sum_{k=1}^n
Y_k(t)f(U_k,W_k)\exp(U_k^T\bolds{\beta}+\eta(W_k))$ and
$s_n[\bolds{\beta},\break\eta](t)=s_n(1;\bolds{\beta},\eta)(t)$.
Define
\begin{eqnarray*}
s[f;\bolds{\beta},\eta](t)&=&E\bigl[Y(t)f(U,W)
\exp\bigl(U^T\bolds{\beta}+\eta(W)\bigr)\bigr]\\
&=&\iint f(u,w)e^{u^T\bolds{\beta}+\eta(w)} q(t,u,w)\,du\,dw.
\end{eqnarray*}
For any functions $f$ and $g$, define
%
\begin{eqnarray}\label{eq:v}
V(f,g)&=&\int_{\mathcal{T}} \biggl\{
\frac{s[fg;\bolds{\beta},\eta_0](t)}{s[\bolds{\beta},\eta_0](t)}
-\frac{s[f;\bolds{\beta},\eta_0](t)}{s[\bolds{\beta},\eta_0](t)}
\frac{s[g;\bolds{\beta},\eta_0](t)}{s[\bolds{\beta},
\eta_0](t)} \biggr\}\nonumber\\[-8pt]\\[-8pt]
&&\hspace*{10.29pt}{}\times
s[\bolds{\beta},\eta_0](t)\,d\Lambda_0(t). \nonumber
\end{eqnarray}
Write $V(f)\equiv V(f,f)$.
Let $\hat{\eta}$ be the estimate that minimizes (\ref{eq:plk}) in
$\mathcal{H}_n$. Then we have
the following theorem.
\begin{theorem}\label{th:eta}
Under conditions \textup{\hyperlink{con:beta}{A1}--\hyperlink{con:lam}{A7}} in the \hyperref[app]{Appendix},
\[
(V+\lambda J)(\hat{\eta}^*-\eta_0)=O_p\bigl(n^{-r/(r+1)}\bigr)
\quad\mbox{and}\quad
(V+\lambda J)(\hat{\eta}-\eta_0)=O_p\bigl(n^{-r/(r+1)}\bigr).
\]
\end{theorem}

This is the optimal convergence rate for estimate of a nonparametric
function. In the view of Lemma \ref{bigsbd}, this theorem also
indicates the same convergence rate in terms of the $L_2$-norm. Also
note that although a higher order of $q_n$ such as $O(n)$ would
yield the same convergence rate for $\hat{\eta}$, it will make the
function space $\mathcal{H}_n$ too big to apply an entropy bound
result that is critical in the proof of Theorem~\ref{th:beta}.

Let
$\mathcal{L}_P(\bolds{\beta})=l_p(\bolds{\beta})
-n\sum_{j=1}^dp_{\theta_j}(|\beta_j|)$,
where
\[
l_p(\bolds{\beta})=\sum_{i=1}^n \int
\{U_i^T\bolds{\beta}+\hat{\eta}(W_i)
-\log s_n[\bolds{\beta},\hat{\eta}](t)\} \,dN_i(t).
\]
Let $\bolds{\beta}_0=(\beta_{10},\ldots,\beta_{d0})^{T}
=(\bolds{\beta}_{10}^{T},\bolds{\beta}_{20}^{T})^T$ be the
true coefficient vector. Without loss of generality, assume that
$\bolds{\beta}_{20}=\mathbf{0}$. Let $s$ be the number of
nonzero components in $\bolds{\beta}_{0}$. Define
$a_n=\max_{j}\{|p'_{\theta_j}(|\beta_{j0}|)|\dvtx\beta_{j0}\ne0\}$,
$b_n=\max_{j}\{|p''_{\theta_j}(|\beta_{j0}|)|\dvtx\break\beta_{j0}\ne0\}$,
and
\begin{eqnarray*}
\mathbf{b}&=& (
p'_{\theta_1}(\beta_{10})\operatorname{sgn}(\beta_{10}),\ldots,
p'_{\theta_s}(\beta_{s0})\operatorname{sgn}(\beta_{s0}) )^{T},\\
\Sigma_{\bolds{\theta}}&=&
\operatorname{diag} \bigl(p'_{\theta_1}(|\beta_{10}|)/|\beta_{10}|,
\ldots,p'_{\theta_s}(|\beta_{s0}|)/|\beta_{s0}| \bigr).
\end{eqnarray*}
Define $\pi_1\dvtx\mathbb{R}^{d+q}\rightarrow\mathbb{R}^{s}$ such that
$\pi_1(u,w)=u_1$, where $u_1$ is the vector of the first $s$
components of $u$. Let $V_0(\pi_1)$ be defined like $V(\pi_1)$ in
(\ref{eq:v}) but with $\bolds{\beta}$ replaced by
$\bolds{\beta}_0$.
\begin{theorem}\label{th:beta}
Under conditions \textup{\hyperlink{con:beta}{A1}--\hyperlink{con:lam}{A7}} in the \hyperref[app]{Appendix},
if $a_n=O(n^{-1/2})$, $b_n=o(1)$ and
$q_n=o(n^{1/2})$, then:
\begin{longlist}
\item
There exists a local maximizer $\hat{\bolds{\beta}}$
of $\mathcal{L}_P(\bolds{\beta})$ such that
$\|\hat{\bolds{\beta}}-\bolds{\beta}_0\|=O_p(n^{-1/2})$.
\item
Further assume that for all $1\le j\le d$,
$\theta_j=o(1)$, $\theta_j^{-1}=o(n^{1/2})$, and
\[
\mathop{\lim\inf}_{n\rightarrow\infty}\,\mathop{\lim\inf}_{u\rightarrow0^{+}}\,
\theta_{j}^{-1}p'_{\theta_{j}}(u)>0.
\]
With probability approaching one, the root-$n$
consistent estimator $\hat{\bolds{\beta}}$ in \textup{(i)} must satisfy
$\hat{\bolds{\beta}}_2=\mathbf{0}$
and
\[
\sqrt{n}\bigl(V_0(\pi_1)+\Sigma_{\bolds{\theta}}\bigr) \bigl\{
\hat{\bolds{\beta}}_1-\bolds{\beta}_{10}+
\bigl(V_0(\pi_1)+\Sigma_{\bolds{\theta}}\bigr)^{-1}\mathbf{b} \bigr\}
\rightarrow N(\mathbf{0},V_0(\pi_1)).
\]
\end{longlist}
\end{theorem}

\subsection{Miscellaneous issues}
\label{ssec:misc}

In this section, we will propose the standard error estimates for
both the parametric and the nonparametric components, and discuss
the selection of the smoothing parameters $\theta$ and $\lambda$.

\subsubsection{Standard error estimates}
Let $l_p(\bolds{\beta})$ be the profile log partial likelihood
in the last iteration of step 3 and
\[
\Sigma_{\bolds{\theta}}(\bolds{\beta})=\operatorname{diag}\{
p'_{\theta_{j}}(|\beta_{1}|)/|\beta_1|,\ldots,
p'_{\theta_{j}}(|\beta_{p}|)/|\beta_p|\}.
\]
Then the standard errors for the nonzero coefficients of
$\hat{\bolds{\beta}}$ are given by the sandwich formula
\[
\hat{\operatorname{cov}}(\hat{\bolds{\beta}})=
\{\nabla^2 l_p(\hat{\bolds{\beta}})-n\Sigma_{\bolds{\theta}}
(\hat{\bolds{\beta}})\}^{-1}
\hat{\operatorname{cov}}\{\nabla l_p(\hat{\bolds{\beta}})\}
\{\nabla^2 l_p(\hat{\bolds{\beta}})-n\Sigma_{\bolds{\theta}}
(\hat{\bolds{\beta}})\}^{-1}.
\]
Sometimes the standard errors for zero coefficients are also of interest.
A discussion of this problem is in Section \ref{sec:disc}.

In (\ref{equ:pplk1}), $\eta$ can be decomposed as
$\eta=\eta^{[0]}+\eta^{[1]}$ where $\eta^{[0]}$ lies in the null
space of the penalty $J$ representing the lower order part and
$\eta^{[1]}$ lies in the complement space representing the higher
order part. A Bayes model interprets (\ref{equ:pplk1}) as a
posterior likelihood when $\eta^{[0]}$ is assigned an improper
constant prior and $\eta^{[1]}$ is assigned a Gaussian prior with
zero mean and certain covariance matrix. The minimizer $\hat{\eta}$
of (\ref{equ:pplk1}) then becomes the posterior mode. When the
minimization is carried out in a data-adaptive function space
$\mathcal{H}_n$ with basis functions
$\bolds{\psi}=(\psi_1,\ldots,\psi_{q_n})^{T}$, we can write
$\hat{\eta}=\bolds{\psi}^{T}\hat{\mathbf{c}}$. Then a
quadratic approximation to (\ref{equ:pplk1}) yields an approximate
posterior covariance matrix for $\mathbf{c}$, which can be used to
construct point-wise confidence intervals for $\eta$.

\subsubsection{Smoothing parameter selection}
As shown in \cite{zou07}, the effective degrees of freedom
for $l_{1}$-penalty model is well approximated by the number of
nonzero coefficients. Note that our SCAD procedure is implemented
by a LASSO approximation at each step. Hence, if we let
$\hat{\mathcal{A}}$ be the set of nonzero coefficients, the AIC
score for selecting $\theta$ in step 3 is
\[
\mathrm{AIC}=2l_{p}(\hat{\bolds{\beta}})+2|\hat{\mathcal{A}}|,
\]
where $|\hat{\mathcal{A}}|$ is the cardinality of $\hat{\mathcal{A}}$.

As illustrated in Section \ref{ssec:npc}, the estimation of $\eta$ in
step 2 can be cast as a density estimation problem with biased sampling.
Let $\operatorname{KL}(\eta_{0},\hat{\eta}_{\lambda})$ be the Kullback--Leibler
distance, as defined in (\ref{equ:ekl}),
between the true ``density'' $e^{\eta_{0}}/\int e^{\eta_{0}}\,d\rmP_{n}^{w}$
and the estimate
$e^{\hat{\eta}_{\lambda}}/\int e^{\hat{\eta}_{\lambda}}\,d\rmP_{n}^{w}$.
An optimal $\lambda$ should minimize
$\operatorname{KL}(\eta_{0},\hat{\eta}_{\lambda})$ or the relative
Kullback--Leibler distance
%
\begin{eqnarray}\label{equ:rkl}
\operatorname{RKL}(\eta_{0},\hat{\eta}_{\lambda})&=&\frac{1}{N}\sum
_{p=1}^{N} \biggl\{
\frac{\int(\eta_{0}(w)-
\hat{\eta}_{\lambda}(w))a_{p}(w)e^{\eta_{0}(w)}\,d\rmP_{n}^{w}}
{\int a_{p}(w)e^{\eta_{0}(w)}\,d\rmP_{n}^{w}}\nonumber\\[-8pt]\\[-8pt]
&&\hspace*{79.3pt}{}+\log\int{a}_{p}(w)e^{\hat{\eta}_{\lambda}(w)}\,d\rmP_{n}^{w}
\biggr\}.\nonumber
\end{eqnarray}
The second term of (\ref{equ:rkl}) is directly computable from the
estimate $\hat{\eta}_{\lambda}$. But the first term needs to be
estimated. Let $\bolds{\psi}$ be the vector of spline basis
functions as in the previous subsection and
$\eta=\bolds{\psi}^{T}\mathbf{c}$. Through a delete-one
cross-validation approximation, a proxy for (\ref{equ:rkl}) can be
derived as
\[
-\frac{1}{N}\sum_{p=1}^{N} \biggl\{
\eta(W_{i})-\log\int{a}_{p}(w)e^{\eta(w)}\,d\rmP_{n}^{w} \biggr\}
+\frac{\operatorname{tr}(P_{\mathbf{1}}Q^{T}H^{-1}QP_{\mathbf{1}})}{N(N-1)},
\]
where $P_{\mathbf{1}}=I-\mathbf{1}\mathbf{1}^{T}/N$,
$Q=(\bolds{\psi}(W_{i_1}),\ldots,\bolds{\psi}(W_{i_{N}}))$,
and $H$ is the hessian matrix for minimizing (\ref{equ:pplk1}) with
respect to the coefficient vector $\mathbf{c}$. $\lambda$ is chosen
to minimize this score.

\section{Numerical studies}
\label{sec:num}

In the simulations, we generated failure times from the exponential
hazard model with
$h(t|U,W)=\exp[U^{T}\bolds{\beta}_{0}+\eta_{0}(W)]. $ We used
the same settings for the parametric component, which consists of
eight covariates $U_{j}, j=1,\ldots,8$. The $U_{j}$'s were
generated from a multivariate normal distribution with zero mean and
$\operatorname{Cov}(U_{j},U_{k})=0.5^{|j-k|}$. The true coefficient vector
was $\bolds{\beta}_{0}=(0.8,0,0,1,0,0,0.6,0)^{T}$.

The theory in Section \ref{ssec:asym} gives the sufficient order for
$q_n$, the number of knots in our smoothing spline estimation of
$\eta$. In practice, \cite{kg04} suggested $q_n=kn^{2/(2m+1)}$
with $k=10$ if the tensor product splines of order $m$ are used.
Since we use tensor product cubic splines in all the simulations
below, our choice is $q_n=10n^{2/5}$.

\subsection{Variable selection for parametric components}
\label{ssec:vsp}

The nonparametric part had one covariate $W$ generated from
$\operatorname{Uniform}(0,1)$. Two different $\eta_{0}$ were used:
\[
\eta_{0a}(w) = 1.5\sin\biggl(2\pi{w}-\frac{\pi}{2}\biggr)
\quad\mbox{or}\quad
\eta_{0b}(w)= 4(w-0.3)^{2}+4.7e^{-w}-3.4643.
\]
Note that both functions satisfies $\int_{0}^{1}\eta_{0}(w)\,dw=0$.
Given $U$ and $W$, the censoring times were generated from
exponential distributions such that the censoring rates are,
respectively, $23\%$ and $40\%$. Sample sizes $n=150$ and 500 were
considered. One thousand data replicates were generated for each of
the four combinations of $\eta_{0}$ and $n$.

For a prediction procedure $\mathcal{M}$ and the estimator
$(\hat{\bolds{\beta}}_{\mathcal{M}},\hat{\eta}_{\mathcal{M}})$
yielded from the procedure,
an appropriate measure for the goodness-of-fit under Cox model with
$h_{0}(t)\equiv{1}$ is the model error:
$\operatorname{ME}(\hat{\bolds{\beta}}_{\mathcal{M}},\hat{\eta
}_{\mathcal{M}})
=E [(\exp(-U^{T}\hat{\bolds{\beta}}_{\mathcal{M}}-\hat{\eta}_{\mathcal{M}}(W))
-\exp(-U^{T}\bolds{\beta}_{0}-\eta_{0}(W)))^{2} ].$ The
relative model error (RME) of $\mathcal{M}_{1}$ versus
$\mathcal{M}_{2}$ is defined as the ratio
$\operatorname{ME}(\hat{\bolds{\beta}}_{\mathcal{M}_{1}},\hat{\eta
}_{\mathcal{M}_{1}})
/\operatorname{ME}(\hat{\bolds{\beta}}_{\mathcal{M}_{2}},\hat{\eta
}_{\mathcal{M}_{2}})$.
The procedure $\mathcal{M}_0$
with complete oracle is used as our benchmark. In $\mathcal{M}_0$,
$(U_{1},U_{4},U_{7}, W)$ are known to be the only contributing covariates,
the exact form of $\eta_{0}$ is known, and the only parameters to be estimated
are the coefficients of $U_1,U_4,U_7$. Note that $\mathcal{M}_0$ can be
implemented only in simulations, but is unrealistic in practice since
neither the contributing covariates nor the form of $\eta_0$ would be known.
We then compare
the performance of the following four procedures,
including the proposed procedures, through their RMEs versus $\mathcal{M}_0$:
\begin{enumerate}[$\mathcal{M}_D$:]
\item[$\mathcal{M}_A$:]
procedure with partial oracle and misspecified parametric $\eta_{0}$,
that is, $(U_{1},U_{4},U_{7},W)$ are known
to be the only contributing covariates but $\eta_{0}$ is misspecified
to be of the parametric form $\eta_{0}(W)=\beta_{W}W$ and $\beta_{W}$
is estimated together with the coefficients for $(U_{1},U_{4},U_{7})$;
\item[$\mathcal{M}_B$:]
procedure with partial oracle and estimated $\eta_{0}$, that is,
$(U_{1},U_{4},U_{7},W)$ are known
to be the only contributing covariates but the form of $\eta_{0}$ is unknown,
and $\eta_0$ is estimated together with the coefficients for
$(U_{1},U_{4},U_{7})$
by penalized profile partial likelihood;
\item[$\mathcal{M}_C$:]
the proposed partial linear procedure with the SCAD penalty on $\bolds
{\beta}$;
\item[$\mathcal{M}_D$:]
the proposed partial linear procedure with the adaptive LASSO penalty
on~$\bolds{\beta}$.
\end{enumerate}
Procedure $\mathcal{M}_A$ has a misspecified covariate effect. We
intend to
show that the estimation results can be unsatisfactory if the semiparametric
form of covariate effect is mistakenly specified as parametric. Procedure
$\mathcal{M}_B$ is ``partial oracle'' and expected to have equal or better
performance than procedures $\mathcal{M}_C$ and $\mathcal{M}_D$. Note, however,
$\mathcal{M}_B$ is unrealistic in practice since the contributing
covariates would not be known. $\mathcal{M}_C$ and $\mathcal{M}_D$
are two versions of the proposed partial linear
procedure with different penalties on
$\bolds{\beta}$.

For each combination of $\eta_{0}$ and $n$, we computed the
following quantities out of the 1000 data replicates: the median
RMEs of the complete oracle
procedure $\mathcal{M}_0$ versus the procedures $\mathcal{M}_A$ to
$\mathcal{M}_D$,
the average number of correctly
selected nonzero coefficients (CC), the average number of incorrectly
selected nonzero coefficients (IC), the proportion of under-fit
replicates that excluded any nonzero coefficients, the proportion of
correct-fit replicates that selected the exact subset model, and the
proportion over-fit replicates that included all three significant
variables and some noise variables. The results are summarized in
Table \ref{tab:vsp}. In general, a partial oracle with misspecified
parametric $\eta_{0}$ (procedure $\mathcal{M}_A$) has much inferior performance
when comparing with the other three procedures; the proposed procedure
with the SCAD penalty (procedure $\mathcal{M}_C$) or the adaptive LASSO penalty
(procedure $\mathcal{M}_D$) is competitive to the partial oracle with estimated
$\eta_{0}$ (procedure $\mathcal{M}_B$); the SCAD penalty performs
slightly better than the
adaptive LASSO penalty. Also, the proposed procedure generally
performs as well as the complete oracle. For procedure $\mathcal{M}_C$,
we also did
some extra computation to evaluate the proposed standard error
estimate of $\bolds{\beta}$. In Table \ref{tab:se}, $\mathrm{SD}$ is
the median absolute deviation divided by 0.6745 of the 1000
nonzero $\hat{\bolds{\beta}}$'s that can be regarded as the true
standard error, $\mathrm{SD}_{m}$ is the median of the 1000 estimated SDs,
and $\mathrm{SD}_{\mathrm{mad}}$ is the median absolute deviation of the 1000
estimated SDs divided by 0.6745. The standard errors were set to 0
for the coefficients estimated as 0s. The results in
Table \ref{tab:se} suggests a good performance of the proposed
standard error formula for $\bolds{\beta}$.

%
\begin{table}
\caption{Variable selection for parametric components (Section \protect
\ref{ssec:vsp})}\label{tab:vsp}
\begin{tabular*}{\tablewidth}{@{\extracolsep{\fill}}l c c c c c c c@{}}
\hline
& & \multicolumn{2}{c}{\textbf{No. of nonzeros}} & &
\multicolumn{3}{c@{}}{\textbf{Proportion of}}\\
[-4pt]
& & \multicolumn{2}{c}{\hrulefill} & &
\multicolumn{3}{c@{}}{\hrulefill}\\
\textbf{Procedure} & \multicolumn{1}{c}{\textbf{MRME}} & \multicolumn{1}{c}{\textbf{CC}}
& \multicolumn{1}{c}{\textbf{IC}} & &
\multicolumn{1}{c}{\textbf{Under-fit}}
& \multicolumn{1}{c}{\textbf{Correct-fit}}
& \multicolumn{1}{c@{}}{\textbf{Over-fit}}\\
\hline
& \multicolumn{7}{c}{$n=150$, $\eta_{0}=\eta_{0a}$ (23\% censoring)}\\
[4pt]
$\mathcal{M}_A$ & 0.168 & -- & -- & & -- & -- & --\\
$\mathcal{M}_B$ & 0.475 & -- & -- & & -- & -- & --\\
$\mathcal{M}_C$ & 0.409 & 2.998 & 0.825 & & 0.002 & 0.476 & 0.522\\
$\mathcal{M}_D$ & 0.387 & 2.998 & 0.959 & & 0.002 & 0.444 & 0.554\\
[4pt]
& \multicolumn{7}{c}{$n=150$, $\eta_{0}=\eta_{0b}$ (40\% censoring)}\\
[4pt]
$\mathcal{M}_A$ & 0.167 & -- & -- & & -- & -- & --\\
$\mathcal{M}_B$ & 0.711 & -- & -- & & -- & -- & --\\
$\mathcal{M}_C$ & 0.518 & 2.996 & 0.949 & & 0.004 & 0.430 & 0.566\\
$\mathcal{M}_D$ & 0.563 & 2.998 & 1.131 & & 0.002 & 0.378 & 0.620\\
[4pt]
& \multicolumn{7}{c}{$n=500$, $\eta_{0}=\eta_{0a}$ (23\% censoring)}\\
[4pt]
$\mathcal{M}_A$ & 0.056 & -- & -- & & -- & -- & --\\
$\mathcal{M}_B$ & 0.431 & -- & -- & & -- & -- & --\\
$\mathcal{M}_C$ & 0.396 & 3.000 & 0.717 & & 0.000 & 0.525 & 0.475\\
$\mathcal{M}_D$ & 0.375 & 3.000 & 0.736 & & 0.000 & 0.540 & 0.460\\
[4pt]
& \multicolumn{7}{c}{$n=500$, $\eta_{0}=\eta_{0b}$ (40\% censoring)}\\
[4pt]
$\mathcal{M}_A$ & 0.057 & -- & -- & & -- & -- & --\\
$\mathcal{M}_B$ & 0.712 & -- & -- & & -- & -- & --\\
$\mathcal{M}_C$ & 0.619 & 3.000 & 0.749 & & 0.000 & 0.512 & 0.488\\
$\mathcal{M}_D$ & 0.628 & 3.000 & 0.776 & & 0.000 & 0.529 & 0.471\\
\hline
\end{tabular*}
\end{table}
%

%
\begin{table}[b]
\caption{Standard deviations for $\hat{\bolds{\beta}}$'s in the
partial linear SCAD procedure $\mathcal{M}_C$ (Section \protect\ref{ssec:vsp})}
\label{tab:se}
\begin{tabular*}{\tablewidth}{@{\extracolsep{\fill}}l  c c  c c  c c@{}}
\hline
& \multicolumn{2}{c}{$\bolds{\hat{\beta}_1}$} &
\multicolumn{2}{c}{$\bolds{\hat{\beta}_4}$}
& \multicolumn{2}{c@{}}{$\bolds{\hat{\beta}_7}$}\\[-4pt]
& \multicolumn{2}{c}{\hrulefill} &
\multicolumn{2}{c}{\hrulefill}
& \multicolumn{2}{c@{}}{\hrulefill}\\
\textbf{\textit{n}, censor \%} & $\mathbf{SD}$ & $\mathbf{SD}_{\bolds m}\ \bolds{(}\mathbf{SD}_{\mathbf{mad}}\bolds{)}$
& $\mathbf{SD}$ & $\mathbf{SD}_{\bolds m}\ \bolds{(}\mathbf{SD}_{\mathbf{mad}}\bolds{)}$ &
$\mathbf{SD}$ & $\mathbf{SD}_{\bolds m}\ \bolds{(}\mathbf{SD}_{\mathbf{mad}}\bolds{)}$\\
\hline
150, 23\% & 0.124 & 0.113 (0.015) & 0.141 & 0.121 (0.017) & 0.135 &
0.109 (0.015)\\
150, 40\% & 0.159 & 0.135 (0.017) & 0.188 & 0.145 (0.021) & 0.155 &
0.128 (0.019)\\
500, 23\% & 0.065 & 0.059 (0.005) & 0.073 & 0.063 (0.005) & 0.062 &
0.057 (0.005)\\
500, 40\% & 0.075 & 0.070 (0.006) & 0.088 & 0.076 (0.006) & 0.078 &
0.066 (0.006)\\
\hline
\end{tabular*}
\end{table}
%

%
\begin{figure}

\includegraphics{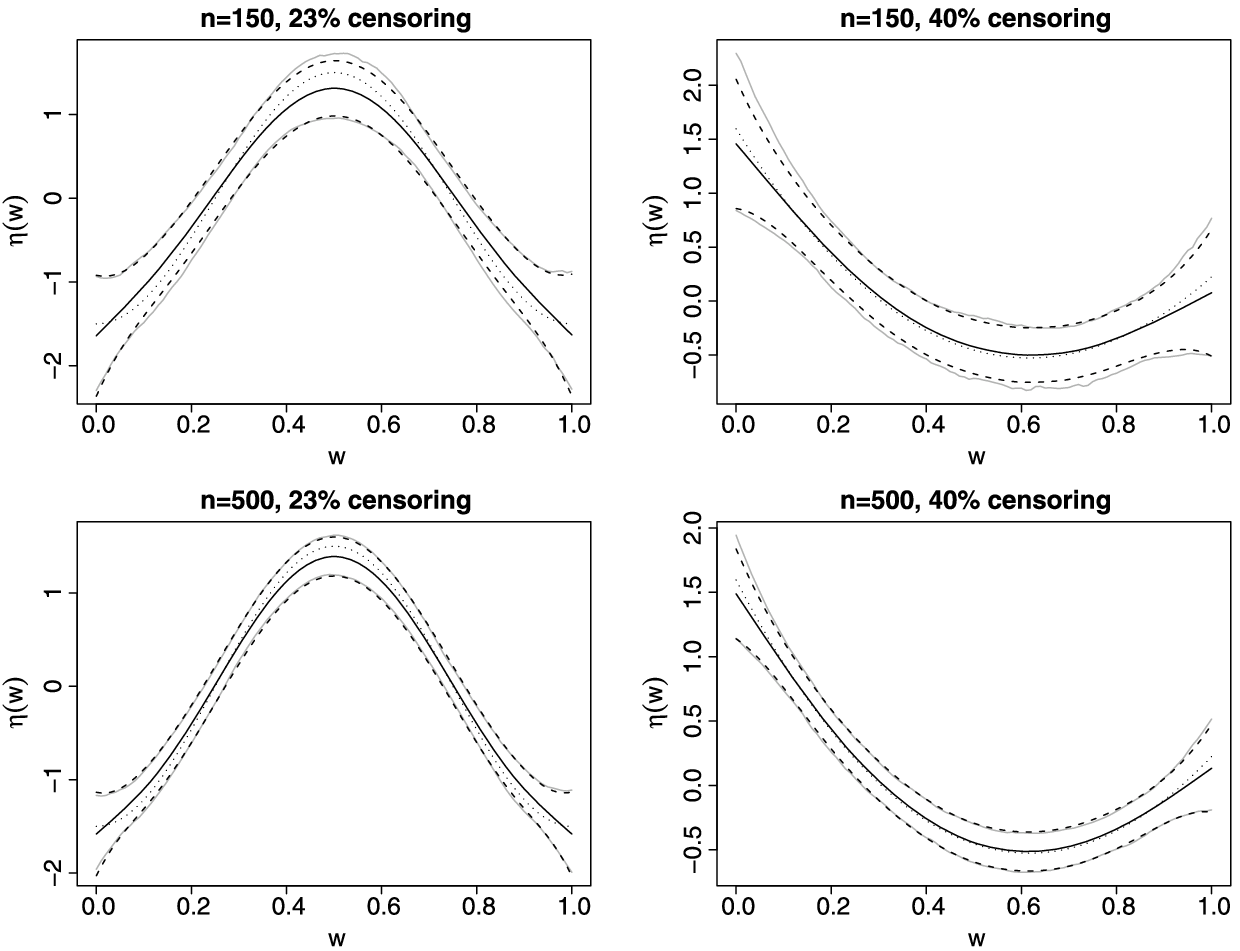}

\caption{Estimates for nonparametric components
(Section \protect\ref{ssec:vsp}).
Dotted lines are true function,
solid lines are connected point-wise mean estimates, faded lines are
connected 0.025 and 0.975 quantiles of the point-wise estimates, and
dashed lines
are the connected point-wise 95\% confidence intervals.}
\label{fig:par}
\end{figure}

To examine the estimation of $\eta_{0}$, we computed the point-wise
estimates at the grid $w=(0,1,\mbox{by}=0.01)$ for each data
replicate. Then at each grid point, the mean, the 0.025 and the
0.975 quantiles of the 1000 estimates, together with the mean of the
1000 95\% confidence intervals were computed. The results are in
Figure \ref{fig:par}. The plots show satisfactory nonparametric
fits and standard error estimates.\looseness=1

\subsection{Model selection for nonparametric components}
\label{ssec:msnp}

In this section, we present some simulations to evaluate the model
selection tool for nonparametric part introduced in
Section \ref{ssec:npc}. We used the SCAD penalty on the parametric
components in this section. Two covariates $W_{1}$ and $W_{2}$,
independently generated from $\operatorname{Uniform}(0,1)$, were used. We
considered two scenarios for the true model of the nonparametric
part: (i) nonparametric univariate model
$\eta_{0}(W)=\eta_{01}(W_1)$ and (ii) nonparametric bivariate
additive model
$\eta_{0}(W)=\eta_{01}(W_{1})+\eta_{02}(W_{2})$. For scenario (i),
the data sets generated in the last section were used, with $W_{1}$
being the existing $W$ covariate and $W_{2}$ being an additional
noise covariate. The fitted model was nonparametric additive in
$W_{1}$ and $W_{2}$. The ratios
$\operatorname{KL}(\hat{\eta},\tilde{\eta})/\operatorname{KL}(\hat{\eta
},\eta_{c})$
for the projections to the univariate models
$\eta_{0}(W)=\eta_{01}(W_1)$ and $\eta_{0}(W)=\eta_{02}(W_2)$ were
computed. For scenario (ii), we considered two sample sizes $n=150$
and 300. The true $\eta_{0}$ was
\[
\eta_{0}(w_{1},w_{2})=0.7\eta_{0a}(w_{1})+0.3\eta_{0b}(w_{2})
\]
or
\[
\eta_{0}(w_{1},w_{2})=\eta_{0a}(w_{1})+\eta_{0b}(w_{2}),
\]
where $\eta_{0a}$ and $\eta_{0b}$ are as defined in
Section \ref{ssec:vsp}. The censoring times were generated from exponential
distributions such that the resulting censoring rates were,
respectively, 25\% and 39\%. Note that both choices of $\eta_{0}$ are
additive in $w_{1}$ and $w_{2}$. The fitted model was the
nonparametric bivariate full model with both the main effects and the
interaction. Then the ratios
$\operatorname{KL}(\hat{\eta},\tilde{\eta})/\operatorname{KL}(\hat{\eta
},\eta_{c})$
for the projections to the bivariate additive model and the two univariate
models were computed. In both scenarios, we claim a reduced model is
feasible when the corresponding ratio
$\operatorname{KL}(\hat{\eta},\tilde{\eta})/\operatorname{KL}(\hat{\eta
},\eta_{c})<0.05$.

%
\begin{table}
\caption{Model selection for nonparametric components (Section \protect
\ref{ssec:msnp})}\label{tab:vsnp}
\begin{tabular*}{\tablewidth}{@{\extracolsep{4in minus 4in}} l c c c c c c c @{}}
\hline
\multirow{2}{30pt}[-6.2pt]{\textbf{Sample size}} &
\multicolumn{3}{c}{\textbf{Proportion of selecting}}
& & \multicolumn{3}{c@{}}{\textbf{Proportion of}}\\[-4pt]
& \multicolumn{3}{c}{\hrulefill}
& & \multicolumn{3}{c@{}}{\hrulefill}\\
& $\bolds{W_{1}}$ & $\bolds{W_{2}}$ & $\bolds{W_{1}\dvtx W_{2}}$
& & \textbf{Under-fit} & \textbf{Correct-fit} & \multicolumn{1}{c@{}}{\textbf{Over-fit}}\\
\hline
& \multicolumn{7}{c}{True model: $\eta_{0}(w_1,w_2)=\eta_{0a}(w_1)$,
23\% censoring}\\
[4pt]
$n=150$ & 1.000 & 0.036 & -- & & 0.000 & 0.964 & 0.036\\
$n=500$ & 1.000 & 0.002 & -- & & 0.000 & 0.998 & 0.002\\
[4pt]
& \multicolumn{7}{c}{True model: $\eta_{0}(w_1,w_2)=\eta_{0b}(w_1)$,
40\% censoring}\\
[4pt]
$n=150$ & 1.000 & 0.304 & -- & & 0.000 & 0.696 & 0.304\\
$n=500$ & 1.000 & 0.062 & -- & & 0.000 & 0.938 & 0.062\\
[4pt]
& \multicolumn{7}{c}{True model: $\eta_{0}(w_1,w_2)=0.7\eta
_{0a}(w_1)+0.3\eta_{0b}(w_2)$, 25\% censoring}\\
[4pt]
$n=150$ & 1.000 & 0.998 & 0.084 & & 0.002 & 0.914 & 0.084\\
$n=300$ & 1.000 & 1.000 & 0.013 & & 0.000 & 0.987 & 0.013\\
[4pt]
& \multicolumn{7}{c}{True model: $\eta_{0}(w_1,w_2)=\eta_{0a}(w_1)+\eta
_{0b}(w_2)$, 39\% censoring}\\
[4pt]
$n=150$ & 1.000 & 0.672 & 0.201 & & 0.328 & 0.471 & 0.201\\
$n=300$ & 1.000 & 0.616 & 0.096 & & 0.384 & 0.520 & 0.096\\
\hline
\end{tabular*}
\end{table}

For each of the eight combinations of $\eta_{0}$ and $n$, we
simulated 1000 data replicates and computed the proportions of
replicates that produced the following results in the reduced model:
selected the main effect of $W_{1}$, selected the main effect of
$W_{2}$, selected the interaction $W_{1}:W_{2}$, under-fitted the
model by excluding at least one truly significant effect, correctly
fitted the model by reducing to the exact subset model, and
over-fitted the model by including all the truly significant effects
and some irrelevant effects. These proportion results are summarized
in Table \ref{tab:vsnp}. It shows that the variable selection tool
for the nonparametric component works very well. The better
performance appears to be associated with bigger sample sizes and
lower censoring rates.

\section{Example}\label{sec:app}

An example in \cite{klein97} is a study on two sexually
transmitted diseases: gonorrhea and chlamydia. The purpose of the
study was to identify factors that are related to time until
reinfection by gonorrhea or chlamydia given an initial infection of
either disease. A sample of 877 individuals with an initial
diagnosis of gonorrhea or chlamydia were followed for reinfection.
Recorded for each individual were follow-up time, indicator of
reinfection, demographic variables including race (white or black,
$U_{1}$), marital status (divorced/separated, married or single,
$U_{2}$~and~$U_{3}$), age at initial infection ($W_{1}$), years of
schooling ($W_{2}$) and type of initial infection (gonorrhea,
chlamydia or both, $U_4$ and $U_5$), behavior factors at the
initial diagnosis including number of partners in the last 30 days
($U_{6}$), indicators of oral sex within past 12 months and within
past 30 days ($U_7$ and $U_8$), indicators of rectal sex within past
12 months and within past 30 days ($U_9$ and $U_{10}$) and condom
use (always, sometimes or never, $U_{11}$ and $U_{12}$), symptom
variables at time of initial infection including presence of
abdominal pain ($U_{13}$), sign of discharge ($U_{14}$), sign of
dysuria ($U_{15}$), sign of itch ($U_{16}$), sign of lesion
($U_{17}$), sign of rash ($U_{18}$) and sign of lymph involvement
($U_{19}$) and symptom variables at time of examination including
involvement vagina at exam ($U_{20}$), discharge at exam ($U_{21}$)
and abnormal node at exam ($U_{22}$).

\begin{figure}[b]

\includegraphics{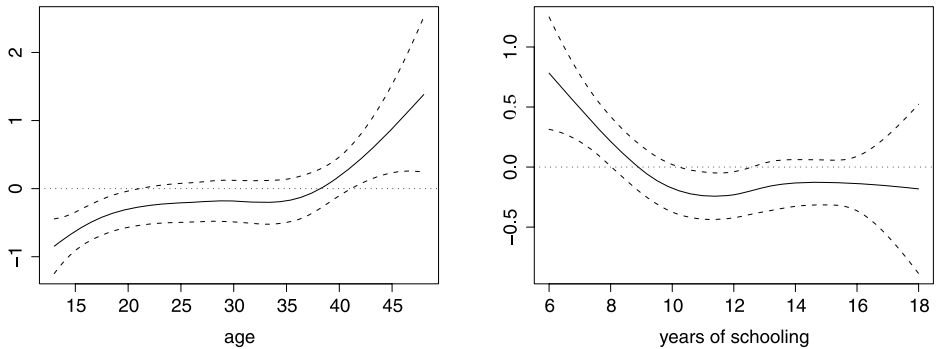}

\caption{Nonparametric component estimates for sexually transmitted
diseases data.
Left: effect of age at initial infection. Right: effect of years of schooling.
Solid lines are the estimates, dashed lines are 95\% confidence intervals
and dotted lines are the reference zero lines.}
\label{fig:std}
\end{figure}

We used $q_n=10\cdot877^{2/5}=151$ knots in all the analysis below.
We first considered the partial linear Cox model
\[
h_{i}(t|Z) = h_{0}(t)\exp\Biggl\{\sum_{j=1}^{3}\eta_{j}(W_{ji})
+\sum_{k=1}^{22}U_{ki}\beta_{k} \Biggr\},
\]
where $\eta_{3}(W_{3i})=\eta_3(W_{1i},W_{2i})$ is the interaction
term between $W_1$ and $W_2$. However, the interaction term was
found to be negligible with the ratio
$\operatorname{KL}(\hat{\eta},\tilde{\eta})/\operatorname{KL}(\hat{\eta
},\eta_{c})=0.003$.
Hence, we took out this interaction term and refitted the model. In
this model, neither $W_1$ (age) nor $W_2$ (years of schooling) in
the nonparametric component were found to be negligible, with the
ratios
$\operatorname{KL}(\hat{\eta},\tilde{\eta})/\operatorname{KL}(\hat{\eta
},\eta_{c})$
equal to $0.633$ for removing $W_{1}$ and $0.259$ for removing
$W_{2}$. Their effects are plotted in Figure \ref{fig:std} together
%
with the 95\% point-wise confidence interval. We can see that the
hazard increased with age at both ends of the age domain (between
age 13 and 20, and between age 38 and 48) and stayed flat in the
middle, and that the hazard decreased with years of school from 6
years to 10 years but stayed flat afterwards. The fitted
coefficients from the proposed method with the SCAD penalty are in
Table \ref{tab:std} together with their standard error estimates.
%
%
\begin{table}
\tabcolsep=3.25pt
\caption{Fitted coefficients and their standard errors
for sexually transmitted diseases data. (Models from top to bottom:
semiparametric relative risk with SCAD penalty and with adaptive LASSO penalty,
parametric relative risk with SCAD penalty and with adaptive LASSO penalty)}
\label{tab:std}
\begin{tabular*}{\textwidth}{@{\extracolsep{\fill}}l c c c c c c c@{}}
\hline
\textbf{\textit{age}} & \textbf{\textit{yschool}} & \textbf{npart} & \textbf{raceW}
& \textbf{maritalM} & \textbf{maritalS} \\
\hline
-- (--) & -- (--) & 0 (--) & 0 (--) & 0 (--) & 0.487 (0.212) \\
-- (--) & -- (--) & \phantom{$-$}0.060 (0.048) & $-$0.127 (0.097) & 0 (--) & 0.448 (0.186) \\
0 (--) & $-$0.059 (0.018) & 0 (--) & 0 (--) & 0 (--) & 0.332 (0.213)\\
0 (--) & $-$0.119 (0.031) & \phantom{$-$}0.026 (0.024) & 0 (--) & 0 (--) & 0.210 (0.119)\\
\hline
\textbf{typeC} & \textbf{typeB} & \textbf{oralY} & \textbf{oralM} & \textbf{rectY} & \textbf{rectM}\\
\hline
$-$0.412 (0.149) & $-$0.337 (0.144) & $-$0.336 (0.201) & $-$0.341 (0.235) & 0 (--) &
0 (--)\\
$-$0.349 (0.137) & $-$0.300 (0.130) & $-$0.330 (0.155) & $-$0.318 (0.173) & 0 (--) &
0 (--)\\
$-$0.376 (0.149) & $-$0.249 (0.145) & $-$0.236 (0.202) & $-$0.348 (0.235) & 0 (--) &
0 (--)\\
$-$0.228 (0.096) & $-$0.083 (0.065) & $-$0.110 (0.058) & $-$0.371 (0.117) & 0 (--) &
0 (--)\\
\hline
\textbf{abdom} & \textbf{disc} & \textbf{dysu} & \textbf{condS} & \textbf{condN} & \textbf{itch} \\
\hline
\phantom{$-$}0.253 (0.151) & 0 (--) & \phantom{$-$}0.193 (0.152) & 0 (--) & $-$0.327 (0.114) & 0 (--) \\
\phantom{$-$}0.177 (0.120) & 0 (--) & \phantom{$-$}0.089 (0.074) & \phantom{$-$}0.152 (0.114) & $-$0.291 (0.106) &
0 (--) \\
\phantom{$-$}0.285 (0.148) & 0 (--) & 0 (--) & 0 (--) & $-$0.296 (0.114) & 0 (--)\\
\phantom{$-$}0.184 (0.094) & 0 (--) & 0 (--) & 0 (--) & $-$0.223 (0.092) & 0 (--)\\
\hline
\textbf{lesion} & \textbf{rash} & \textbf{lymph} & \textbf{involve} & \textbf{discE} & \textbf{node}\\
\hline
0 (--) & 0 (--) & 0 (--) & \phantom{$-$}0.423 (0.166) & $-$0.460 (0.220) & 0 (--)\\
0 (--) & 0 (--) & 0 (--) & \phantom{$-$}0.327 (0.159) & $-$0.407 (0.209) & 0 (--)\\
0 (--) & 0 (--) & 0 (--) & \phantom{$-$}0.392 (0.168) & $-$0.443 (0.221) & 0 (--)\\
0 (--) & 0 (--) & 0 (--) & \phantom{$-$}0.289 (0.133) & $-$0.280 (0.163) & 0 (--)\\
\hline
\end{tabular*}
\end{table}
For comparisons, Table \ref{tab:std} also lists the fitted
coefficients and standard errors for three other models, namely the
proposed semiparametric relative risk model with the adaptive LASSO
penalty, and the parametric relative risk models with the SCAD and
the adaptive LASSO penalties. We can see that the SCAD penalty
yielded sparser models than the adaptive LASSO penalty, and that both parametric
models missed the age effect. Common factors identified by all the
four procedures to be associated with reinfection risk are marital
status, type of infection, oral sex behavior, condom use, sign of
abdominal pain, sign of lymph involvement and sign of discharge at
exam.\looseness=1

\section{Discussion}\label{sec:disc}

We have proposed a Cox PH model with semiparametric relative risk.
The nonparametric part of the risk is estimated by smoothing spline
ANOVA model and model selection procedure derived based on a
Kullback--Leibler geometry. The parametric part of the risk is
estimated by penalized profile partial likelihood and variable
selection achieved by choosing a nonconcave penalty. Both
theoretical and numerical studies show promising results for the
proposed method. An important question in using the method
in practice is which covariate effects should be treated as
parametric. We suggest the following guideline for making
choices. As a starting point, the effects of all the continuous covariates
are put in the nonparametric part and those of the discrete covariates
in the parametric part. If the estimation results show that some of
the continuous covariate effects can be described by certain parametric
forms such as linear form, then a new model can be fitted with
those continuous covariate effects moved to the parametric part.
In this way, one can take full advantage of the flexible exploratory
analysis provided by the proposed method.

We thank a referee for raising the interesting question on
the standard error estimates for the coefficients estimated to be 0 in
$\hat{\bolds{\beta}}$. References \cite{lasso96} and \cite{scad01}
suggested to set these standard errors to 0s based on the belief
that those covariates with zero coefficient estimates are not important.
This is the approach adopted here.
When such a belief is in doubt, nonzero standard errors are preferred
even for coefficients estimated to be 0's. This problem has been
addressed only in a few papers.
Reference \cite{osbo00} looked at the problem for LASSO
but it is based on a smooth approximation.
Reference \cite{park08} presented a Bayesian approach and pointed out that
no fully satisfactory frequentist solution had been proposed so far,
no matter LASSO or SCAD variable selection procedure is considered.
This problem presents an interesting challenge that we hope to
address in some future work.

Another choice of $p_{\theta}(|\cdot|)$ is the adaptive LASSO
penalty \cite{zou06}.
Our simulations in
Section \ref{ssec:vsp} indicates a similar performance when compared to
the SCAD penalty. So we decided not to present the details here.

Although our method is presented for time-independent covariates,
a lengthier argument modifying \cite{osu93} can yield similar
theoretical results for external time-dependent covariates
\cite{kal02}. However, the implementation of such
extension is more complicated and not pursued here.

A recently proposed nonparametric component selection procedure in a
penalized likelihood framework is the COSSO method in \cite{cosso06}
where the penalty switches from $J(\eta)$ to $J^{1/2}(\eta)$.
Taking advantage of the smoothing spline ANOVA decomposition,
the COSSO method does model selection by applying a soft thresholding
type operation to the function components. An extension of
COSSO to the Cox proportional hazards model with nonparametric relative
risk is available in \cite{leng06}.
Although a similar extension to our proportional hazards model with
semiparametric relative risk is of interest,
it is not clear whether the theoretical properties of the COSSO method
such as the existence and the convergence rate of the COSSO estimator
can be transferred to the estimation of $\eta$ under
our semiparametric setting.
Furthermore, the dimension of the function space in COSSO is $O(n)$,
too big to allow an entropy bound that is critical in deriving
the asymptotic properties of $\hat{\bolds{\beta}}$.

\begin{appendix}\label{app}
\section*{Appendix: Proofs}

For $z=(u,w)$, let $p_z(t)=\exp(u^T\bolds{\beta}
+\eta_0(w))q(t,u,w)/s[\bolds{\beta},\eta_0](t)$ and
$\bar{f}_t\equiv\iint f(u,w)p_z(t)\,du\,dw$.
Let $\tilde{S}(t,u,w)=E[Y(t)=1|U=u,W=w]=P(Y(t)=1|U=u,W=w)$ and
$q(t,u,w)=\tilde{S}(t,u,w)p(u,w)$, where $p(u,w)$ is the density function
of $(U,W)$. Let $\mathcal{Z}=\mathcal{U}\times\mathcal{W}$
be the domain of the covariate
$Z=(U^T,W^T)^T$.
We need the following conditions.
\begin{enumerate}[A1.]
\item[A1.]\hypertarget{con:beta}
The true coefficient $\beta_0$ is an interior point of a bounded subset of
$\mathbb{R}^d$.
\item[A2.]\hypertarget{con:cov}
The domain $\mathcal{Z}$ of covariate is a compact set in $\mathbb{R}^{d+q}$.
\item[A3.]\hypertarget{con:cen}
Failure time $T$ and censoring time $C$ are conditionally independent
given the covariate $Z$.
\item[A4.]\hypertarget{con:bound}
Assume the observations are in a finite time interval $[0,\tau]$.
Assume that the baseline hazard function
$h_0(t)$ is bounded away from zero and infinity.
\item[A5.]\hypertarget{con:q}
Assume that there exist constants
$k_2>k_1>0$ such that
$k_1<q(t,u,w)<k_2$ and $ |\frac{\partial}{\partial t}q(t,u,w) |<k_2$.
\item[A6.]\hypertarget{con:eta}
Assume the true function $\eta_0\in\mathcal{H}$.
For any $\eta$ in a sufficiently big convex neighborhood $B_0$ of $\eta_0$,
there exist constants $c_1, c_2 >0$ such that
$c_1e^{\eta_0(w)}\le e^{\eta(w)}\le c_2e^{\eta_0(w)}$ for all $w$.
\item[A7.]\hypertarget{con:lam}
The smoothing parameter $\lambda\asymp n^{-r/(r+1)}$.
\end{enumerate}

Condition \hyperlink{con:beta}{A1} requires that $\beta_0$ is not on the
boundary of the parameter space. Condition \hyperlink{con:cov}{A2} is also a common
boundedness assumption on covariate. Condition~\hyperlink{con:cen}{A3}
assumes noninformative censoring. Condition \hyperlink{con:bound}{A4}
is the common boundedness assumption on the baseline hazard.
Condition \hyperlink{con:q}{A5} bounds the joint density of $(T,Z)$ and thus
also the derivatives of the partial likelihood. Condition
\hyperlink{con:eta}{A6}
assumes that $\eta_0$
has proper level of smoothness and integrates to zero.
The neighborhood
$B_0$ in condition \hyperlink{con:eta}{A6} should be big enough to contain
all the estimates of $\eta_0$ considered below. When the members of
$B_0$ are all uniformly bounded, condition \hyperlink{con:eta}{A6} is
automatically satisfied. The order for $\lambda$
in condition \hyperlink{con:lam}{A7} matches that in standard
smooth spline problems.

We first show the equivalence between $V(\cdot)$ and the $L_2$-norm
$\|\cdot\|_2^2$.
\begin{lemma}\label{bigsbd}
Let $f\in\mathcal{H}$.
Then there exist constants $0<c_3\le c_4<\infty$ such that
\[
c_3\|f\|_2^2 \le V(f)\le c_4\|f\|_2^2.
\]
\end{lemma}
\begin{pf}
For $z=(u,w)$, let $p_z(t)=\exp(u^T\bolds{\beta}
+\eta_0(w))q(t,u,w)/s[\bolds{\beta},\eta_0](t)$ and
$\bar{f}_t\equiv\iint f(u,w)p_z(t)\,du\,dw$. Simple
algebraic manipulation yields
\[
V(f)=\int_{\mathcal{T}} \biggl\{
\iint\bigl(f(u,w)-\bar{f}_t\bigr)^2 p_z(t)\,du\,dw \biggr\}
s[\bolds{\beta},\eta_0](t)\,d\Lambda_0(t).
\]
By conditions \hyperlink{con:bound}{A4} and \hyperlink{con:q}{A5}, there exist
positive constants $c_1$ and $c_2$ such that
\begin{eqnarray*}
&&c_1\int_{\mathcal{T}} \biggl\{\iint
\bigl(f(u,w)-\bar{f}_t\bigr)^2\,du\,dw \biggr\}\,d\Lambda_0(t)\\
&&\qquad\le V(f)
\le
c_2\int_{\mathcal{T}} \biggl\{\iint
\bigl(f(u,w)-\bar{f}_t\bigr)^2\,du\,dw \biggr\}\,d\Lambda_0(t).
\end{eqnarray*}
Let $m(\mathcal{Z})<\infty$ be the Lebesgue measure of $\mathcal{Z}$. Then
\begin{eqnarray*}
&&\int_{\mathcal{T}} \biggl\{\iint\bigl(f(u,w)-\bar{f}_t\bigr)^2\,du\,dw \biggr\}\,d\Lambda_0(t)\\
&&\qquad=\Lambda_0(\tau)\iint f^2(u,w)\,du\,dw
+m(\mathcal{Z})\iint[\bar{f}_t]^2\,d\Lambda_0(t).
\end{eqnarray*}
The lemma follows from
the Cauchy--Schwarz inequality and condition \hyperlink{con:bound}{A4}.
\end{pf}
\begin{pf*}{Proof of Theorem \ref{th:eta}} We will prove the results
using an eigenvalue analysis of three steps. In the first step
(linear approximation), we show the convergence rate
$O_p(n^{-r/(r+1)})$ for the minimizer $\tilde{\eta}$ of a quadratic
approximation to~(\ref{eq:plk}). In the second step (approximation
error), we show that the difference between $\tilde{\eta}$ and the
estimate $\hat{\eta}^*$ in $\mathcal{H}$ is also
$O_p(n^{-r/(r+1)})$, and so is the convergence rate of
$\hat{\eta}^*$. In the third step (semiparametric approximation), we
show that the projection $\eta^*$ of $\hat{\eta}^*$ in
$\mathcal{H}_n$ is not so different from either $\hat{\eta}^*$ or
the estimate $\hat{\eta}$ in $\mathcal{H}_n$, and then the
convergence rate of $\hat{\eta}$ follows.

A quadratic function
$B$ is said to be \textit{completely continuous} with respect to
another quadratic functional $A$, if for any $\varepsilon>0$, there
exists a finite number of linear functionals $L_1,\ldots,L_k$ such
that $L_j f=0, j=1,\ldots,k$, implies that $B(f)\le\varepsilon A(f)$;
When a quadratic functional $B$ is \textit{completely continuous}
with respect to another quadratic functional $A$,
there exists eigenfunctions $\{\phi_\nu,\nu=1,2,\ldots\}$
such that $B(\phi_\nu,\phi_\mu)=\delta_{\nu\mu}$ and
$A(\phi_\nu,\phi_\mu)=\rho_\nu\delta_{\nu\mu}$, where
$\delta_{\nu\mu}$ is the Kronecker delta and
$0\le\rho_\nu\uparrow\infty$. And functions satisfying $A(f)<\infty$
can be expressed as a \textit{Fourier series expansion}
$f=\sum_\nu f_\nu\phi_\nu$, where $f_\nu=B(f,\phi_\nu)$ are the
\textit{Fourier coefficients}.
See, for example, \cite{gu02} and
\cite{weinberger74}.

We first present two lemmas without proof. The first one follows
directly from the results in Section 8.1 of \cite{gu02} and
Lemma \ref{bigsbd}. The second one is exactly Lem\-ma~8.1 in
\cite{gu02}.
\begin{lemma}\label{lem:vj}
$V$ is completely continuous to $J$ and the eigenvalues $\rho_\nu$
of $J$ with respect to $V$ satisfy that
as $\nu\rightarrow\infty$,
$\rho_\nu^{-1}=O(\nu^{r})$.
\end{lemma}
\begin{lemma}\label{sumbd}
As $\lambda\rightarrow{0}$, the sums
$\sum_{\nu}\frac{\lambda\rho_{\nu}}{(1+\lambda\rho_{\nu})^{2}}$,
$\sum_{\nu}\frac{1}{(1+\lambda\rho_{\nu})^{2}}$,
and $\sum_{\nu}\frac{1}{1+\lambda\rho_{\nu}}$ are all of
order $O(\lambda^{-1/r})$.
\end{lemma}
\begin{Step}[(Linear approximation)]
A linear approximation $\tilde{\eta}$ to $\hat{\eta}^*$ is the
minimizer of a quadratic approximation to (\ref{eq:plk}),
%
\begin{eqnarray}\label{eq:quad}
&&-\frac{1}{n}\sum_{i=1}^{n}\int_{\mathcal{T}} \{\eta(W_i)-
s_n^{-1}[\bolds{\beta},\eta_0](t)s_n[\eta-\eta_0;
\bolds{\beta},\eta_0](t) \}\,dN_i(t)\nonumber\\[-8pt]\\[-8pt]
&&\qquad{} + \frac{1}{2}V(\eta-\eta_0)+\frac{\lambda}{2}J(\eta).\nonumber
\end{eqnarray}
Let $\eta=\sum_\nu\eta_\nu\phi_\nu$ and
$\eta_0=\sum_\nu\eta_{\nu,0}\phi_\nu$ be the Fourier expansions of
$\eta$ and $\eta_0$. Plugging them into (\ref{eq:quad}) and dropping the
terms not involving $\eta$ yield
%
\begin{equation}\label{eq:quadf}
\sum_{\nu} \biggl\{
-\eta_\nu\gamma_\nu
+\frac{1}{2}(\eta_\nu-\eta_{\nu,0})^2+\frac{\lambda}{2}\rho_{\nu}\eta
_\nu^2\biggr\},
\end{equation}
where $\gamma_\nu=\frac{1}{n}\sum_{i=1}^{n}
\int_{\mathcal{T}}\{\phi_\nu(W_i)-
s_n^{-1}[\bolds{\beta},\eta_0](t)s_n[\phi_\nu;
\bolds{\beta},\eta_0](t)\}\,dN_i(t)$.
The\break Fourier coefficients that minimize (\ref{eq:quadf}) are
$
\tilde{\eta}_\nu=(\gamma_\nu+\eta_{\nu,0})/(1+\lambda\rho_\nu).
$
Note that $\int\phi_\nu(w)\,dw=0$ and
$V(\phi_\nu)=1$. Straightforward calculation gives
$E[\gamma_{\nu}]=0$ and $E[\gamma_{\nu}^{2}]=n^{-1}$. Then
%
\begin{eqnarray}\label{vj0:8}
E[V(\tilde{\eta}-\eta_{0})]&=&\frac{1}{n}
\sum_{\nu}\frac{1}{(1+\lambda\rho_{\nu})^{2}}
+\lambda\sum_{\nu}\frac{\lambda\rho_{\nu}}{(1+\lambda\rho_{\nu})^{2}}
\rho_{\nu}\eta_{\nu,0}^{2},\nonumber\\[-8pt]\\[-8pt]
E[\lambda{J}(\tilde{\eta}-\eta_{0})]
&=&\frac{1}{n}\sum_{\nu}\frac{\lambda\rho_{\nu}}{(1+\lambda\rho_{\nu})^{2}}
+\lambda\sum_{\nu}\frac{(\lambda\rho_{\nu})^{2}}{(1+\lambda\rho_{\nu})^{2}}
\rho_{\nu}\eta_{\nu,0}^{2}.\nonumber
\end{eqnarray}
Combining Lemma \ref{sumbd} and (\ref{vj0:8}), we obtain that
$(V+\lambda{J})(\tilde{\eta}-\eta_{0})
=O_{p}(\lambda+n^{-1}\lambda^{-1/r})$, as $n\rightarrow\infty$
and $\lambda\rightarrow{0}$.
\end{Step}
\begin{Step}[(Approximation error)]
We now investigate the approximation error $\hat{\eta}^*-\tilde{\eta}$ and
prove the convergence rate of $\hat{\eta}^*$. Define $A_{f,g}(\alpha)$
and $B_{f,g}(\alpha)$, respectively, as the resulting functionals from
setting $\eta=f+\alpha{g}$ in (\ref{eq:plk}) and (\ref{eq:quad}).
Differentiating them with respect to $\alpha$ and then setting
$\alpha=0$ yields
%
\begin{eqnarray}\quad
\label{a0_hzd:8}
\dot{A}_{f,g}(0)
&=& -\frac{1}{n}\sum_{i=1}^{n}\int_{\mathcal{T}} \{g(W_i)
-s_n^{-1}[\bolds{\beta},f](t)
s_n[g;\bolds{\beta},f](t) \}\,dN_i(t)\nonumber\\[-8pt]\\[-8pt]
&&{} +\lambda{J}(f,g),\nonumber
\\
\label{b0_hzd:8}
\dot{B}_{f,g}(0)
&=& -\frac{1}{n}\sum_{i=1}^{n}\int_{\mathcal{T}} \{g(W_i)
-s_n^{-1}[\bolds{\beta},\eta_0](t) s_n[g;
\bolds{\beta},\eta_0](t) \}\,dN_i(t)\nonumber\\[-8pt]\\[-8pt]
&&{} + V(f-\eta_{0},g)+\lambda{J}(f,g).\nonumber
\end{eqnarray}
Set $f=\hat{\eta}^*$ and $g=\hat{\eta}^*-\tilde{\eta}$ in
(\ref{a0_hzd:8}), and set $f=\tilde{\eta}$ and
$g=\hat{\eta}^*-\tilde{\eta}$ in
(\ref{b0_hzd:8}). Then subtracting the resulted equations gives
%
\begin{eqnarray}\label{a1b1_hzd:8}
&&\mu_{\hat{\eta}^*}(\hat{\eta}^*-\tilde{\eta})-\mu_{\tilde{\eta}}(\hat
{\eta}^*-\tilde{\eta})
+\lambda{J}(\hat{\eta}^*-\tilde{\eta})\nonumber\\[-8pt]\\[-8pt]
&&\qquad=V(\tilde{\eta}-\eta_{0},\hat{\eta}^*-\tilde{\eta})
+\mu_{\eta_0}(\hat{\eta}^*-\tilde{\eta})-\mu_{\tilde{\eta}}(\hat{\eta
}^*-\tilde{\eta}),\nonumber
\end{eqnarray}
where $\mu_{f}(g)\equiv\frac{1}{n}\sum_{i=1}^{n}\int_{\mathcal{T}}
s_n^{-1}[\bolds{\beta},f](t) s_n[g;\bolds{\beta},f](t)\,dN_i(t)$. Define
\[
S_n[f,g](t)=\frac{s_n[fg;\bolds{\beta},\eta_0](t)}{s_n[\bolds{\beta
},\eta_0](t)}-
\frac{s_n[f;\bolds{\beta},\eta_0](t)}{s_n[\bolds{\beta},\eta_0](t)}
\frac{s_n[g;\bolds{\beta},\eta_0](t)}{s_n[\bolds{\beta},\eta_0](t)}
\]
and $S[f,g](t)$ be its limit. The following lemma is needed to proceed.
\begin{lemma}\label{vv_hzd:t}
\begin{eqnarray*}
&&\frac{1}{n}\sum_{i=1}^{n}
\int_{\mathcal{T}}S_n[f,g](t)\,dN_i(t)\\
&&\qquad
=V(f,g)+o_{p}\bigl(\{(V+\lambda{J})(f)(V+\lambda{J})(g)\}^{1/2}\bigr).
\end{eqnarray*}
\end{lemma}
\begin{pf}
Let $f=\sum_{\nu}f_{\nu}\phi_{\nu}$ and $g=\sum_{\mu}g_{\mu}\phi_{\mu}$
be the Fourier series expansion
of $f$ and $g$.
Reference \cite{andegill82} shows that 
${\sup_t}|s_n[\bolds{\beta},\eta_0]-s[\bolds{\beta},\eta_0]|$ converges
to zero
in probability.
Note that
\[
M(t)\equiv M(t|Z) =
N(t)-\int_0^t s[\bolds{\beta},\eta_0](\tau)\,d\Lambda_0(\tau)
\]
defines a local martingale with mean zero.
Combining the above uniform convergence result and the martingale
property with the boundedness condition, we obtain that for any $\nu$
and $\mu$,
\[
E \biggl[ \biggl\{\int_{\mathcal{T}}S[\phi_\nu,\phi_\mu](t)\,dN(t)
-V(\phi_{\nu},\phi_{\mu}) \biggr\}^{2} \biggr]<\infty.
\]
Then from the Cauchy--Schwarz inequality and Lemma \ref{sumbd},
\begin{eqnarray*}
&&\Biggl|\frac{1}{n}\sum_{i=1}^{n}\int_{\mathcal{T}}S_n[f,g](t)
\,dN_i(t)-V(f,g) \Biggr|\\
&&\qquad= \Biggl|\sum_{\nu}\sum_{\mu}f_{\nu}g_{\mu}
\Biggl\{\frac{1}{n}\sum_{i=1}^{n}\int_{\mathcal{T}}
S_n[\phi_\nu,\phi_\mu](t)\,dN_i(t)
-V(\phi_{\nu}, \phi_{\mu})\Biggr\} \Biggr|\\
&&\qquad\leq\Biggl\{\sum_{\nu}\sum_{\mu}\frac{1}{1+\lambda\rho_{\nu}}
\frac{1}{1+\lambda\rho_{\mu}}\\
&&\qquad\quad\hspace*{35.3pt}{} \times\Biggl\{\frac{1}{n}\sum_{i=1}^{n}\int_{\mathcal{T}}
S_n[\phi_\nu,\phi_\mu](t)\,dN_i(t)
-V(\phi_{\nu},\phi_{\mu}) \Biggr\}^{2} \Biggr\}^{1/2}\\
&&\qquad\quad{} \times\biggl\{\sum_{\nu}\sum_{\mu}(1+\lambda\rho_{\nu})
(1+\lambda\rho_{\mu})f_{\nu}^{2}g_{\mu}^{2} \biggr\}^{1/2}\\
&&\qquad=O_{p}(n^{-1/2}\lambda^{-1/r})\{(V+\lambda{J})(f)(V+\lambda{J})(g)\}
^{1/2}.
\end{eqnarray*}
\upqed\end{pf}

A Taylor expansion at $\eta_0$ gives
\begin{eqnarray*}
\mu_{\hat{\eta}^*}(\hat{\eta}^*-\tilde{\eta})-\mu_{\tilde{\eta}}(\hat
{\eta}^*-\tilde{\eta})
&=&\frac{1}{n}\sum_{i=1}^{n}\int_{\mathcal{T}}
S_n[\hat{\eta}^*-\tilde{\eta},\hat{\eta}^*-\tilde{\eta
}](t)\,dN_i(t)\bigl(1+o_p(1)\bigr),\\
\mu_{\tilde{\eta}}(\hat{\eta}^*-\tilde{\eta})-\mu_{\eta_0}(\hat{\eta
}^*-\tilde{\eta})
&=&\frac{1}{n}\sum_{i=1}^{n}\int_{\mathcal{T}}
S_n[\tilde{\eta}-\eta_0,\hat{\eta}^*-\tilde{\eta}](t)\,dN_i(t)\bigl(1+o_p(1)\bigr).
\end{eqnarray*}
Then by the mean value theorem, condition \hyperlink{con:eta}{A6},
Lemma \ref{vv_hzd:t} and (\ref{a1b1_hzd:8}),
\begin{eqnarray*}
&&(c_{1}V+\lambda{J})(\hat{\eta}^*-\tilde{\eta})\bigl(1+o_{p}(1)\bigr)\\
&&\qquad\leq\{(|1-c|V+\lambda{J})(\hat{\eta}^*-\tilde{\eta})\}^{1/2}
O_{p}\bigl(\{(|1-c|V+\lambda{J})(\tilde{\eta}-\eta_{0})\}^{1/2}\bigr)
\end{eqnarray*}
for some $c\in[c_{1},c_{2}]$. Then the convergence rate of
$\hat{\eta}^*$ follows from that of
$\tilde{\eta}$ proved in the previous step.
\end{Step}
\begin{Step}[(Semiparametric approximation)]
Our last goal is the convergence rate for the minimizer $\hat{\eta}$
in the space $\mathcal{H}_{n}$.
For any $h\in\mathcal{H}\ominus\mathcal{H}_{n}$,\break one has
$h(W_{i_l})=J(R_{J}(W_{i_l},\cdot),h)=0$,
so $s_{q_n}[h^j;\bolds{\beta},\eta_0](t)=\break\frac{1}{q_n}\sum
_{k=1}^{q_n}Y_{i_k}(t)h^j(W_{i_k})
\exp(U_{i_k}^T\bolds{\beta}+\eta(W_{i_k}))=0$ for $j=1,2$
and\break $\sum_{l=1}^{q_n}\int_{\mathcal{T}}S_{q_n}[h,h](t)\,dN_{i_l}(t)=0$.
Hence, by the same arguments used in
the proof of Lemma \ref{vv_hzd:t},
%
\begin{eqnarray}\label{key_hzd:t}
V(h) &=& \Biggl|\frac{1}{q_n}\sum_{l=1}^{q_n}\int_{\mathcal{T}}
S_{q_n}[h,h](t)\,dN_{i_l}(t)-V(h) \Biggr|\nonumber\\[-8pt]\\[-8pt]
&=&
O_{p}(q_n^{-1/2}\lambda^{-1/r})(V+\lambda{J})(h)=o_{p}(\lambda{J}(h)),\nonumber
\end{eqnarray}
where the last equality follows from $q_n\asymp n^{2/(r+1)+\varepsilon}$
and condition \hyperlink{con:lam}{A7}.

Let $\eta^{*}$ be the projection of $\hat{\eta}^{*}$ in $\mathcal{H}_{n}$.
Setting $f=\hat{\eta}^*$ and $g=\hat{\eta}^*-\eta^{*}$ in
(\ref{a0_hzd:8}) and noting that $J(\eta^{*},\hat{\eta}^*-\eta^{*})=0$,
some algebra yields
%
\begin{eqnarray}\label{gn1_hzd:8}
\lambda{J}(\hat{\eta}^*-\eta^{*})
&=& \Biggl\{\frac{1}{n}\sum_{i=1}^{n}\int_{\mathcal{T}}
(\hat{\eta}^*-\eta^{*})(W_i)\,dN_i(t)-\mu_{\eta_0}(\hat{\eta}^*-\eta^{*})
\Biggr\}\nonumber\\[-8pt]\\[-8pt]
&&{} - \{\mu_{\hat{\eta}^*}(\hat{\eta}^*-\eta^{*})-\mu_{\eta_0}(\hat{\eta
}^*-\eta^{*}) \}.\nonumber
\end{eqnarray}
Recall that
$\gamma_\nu=\frac{1}{n}\sum_{i=1}^{n}\int_{\mathcal{T}}\phi_\nu(W_i)\,dN_i(t)
-\mu_{\eta_0}(\phi_\nu)$
with $E[\gamma_\nu]=0$ and\break $E[\gamma_\nu^2]=1/n$.
An application of the Cauchy--Schwarz inequality\break and Lemma \ref{sumbd}
shows that the first term in (\ref{gn1_hzd:8}) is of order
$\{(V+\lambda{J})(\hat{\eta}^*-\eta^{*})\}^{1/2}O_{p}(n^{-1/2}\lambda^{-1/2r})$.
By the mean value theorem, condition \hyperlink{con:eta}{A6}, Lemma \ref{vv_hzd:t}
and (\ref{key_hzd:t}), the remaining term in (\ref{gn1_hzd:8}) is of order
$o_{p}(\{\lambda{J}(\hat{\eta}^*-\eta^{*})
(V+\lambda{J})(\hat{\eta}^*-\eta_{0})\}^{1/2})$. These, combined with
(\ref{gn1_hzd:8}) and the convergence rates of $\hat{\eta}^*$, yield
$\lambda{J}(\hat{\eta}^*-\eta^{*})=O_{p}(n^{-1}\lambda^{-1/r}+\lambda)$ and
$V(\hat{\eta}^*-\eta^{*})=o_{p}(n^{-1}\lambda^{-1/r}+\lambda)$.

Note that $J(\hat{\eta}^*-\eta^{*},\eta^{*})=J(\hat{\eta}^*-\eta
^{*},\hat{\eta})=0$,
so $J(\hat{\eta}^*,\hat{\eta}^*-\hat{\eta})=J(\hat{\eta}^*-\eta^{*})
+J(\eta^{*},\eta^{*}-\hat{\eta})$.
Set $f=\hat{\eta}$ and $g=\hat{\eta}-\eta^{*}$ in
(\ref{a0_hzd:8}), and set $f=\hat{\eta}^*$ and $g=\hat{\eta}^*-\hat{\eta
}$ in
(\ref{a0_hzd:8}). Adding the resulted equations yields
%
\begin{eqnarray}\label{gngn_hzd:8}
&&\mu_{\hat{\eta}}(\hat{\eta}-\eta^{*})
-\mu_{\eta_0}(\hat{\eta}-\eta^{*})
+\lambda{J}(\hat{\eta}-\eta^{*})
+\lambda{J}(\hat{\eta}^*-\eta^{*})\nonumber\\
&&\qquad= \Biggl\{\frac{1}{n}\sum_{i=1}^{n}\int_{\mathcal{T}}
(\hat{\eta}^*-\eta^{*})(W_i)\,dN_i(t)-\mu_{\eta_0}(\hat{\eta}^*-\eta^{*})
\Biggr\}\nonumber\\[-8pt]\\[-8pt]
&&\qquad\quad{}- \{\mu_{\hat{\eta}^*}(\hat{\eta}^*-\eta^{*})-\mu_{\eta_0}(\hat{\eta
}^*-\eta^{*}) \}\nonumber\\
&&\qquad\quad{}
+ \{\mu_{\hat{\eta}^*}(\hat{\eta}-\eta^{*})-
\mu_{\eta_0}(\hat{\eta}-\eta^{*}) \}.\nonumber
\end{eqnarray}

By the mean value theorem, condition \hyperlink{con:eta}{A6}, and
Lemma \ref{vv_hzd:t}, the left-hand side of (\ref{gngn_hzd:8}) is
bounded from below by
\[
(c_{1}V+\lambda{J})(\hat{\eta}-\eta^{*})\bigl(1+o_{p}(1)\bigr)
+\lambda{J}(\hat{\eta}^*-\eta^{*}).
\]
For the right-hand side, the terms in the first and second brackets
are, respectively, of the orders
$\{(V+\lambda{J})(\hat{\eta}^*-\eta^{*})\}^{1/2}O_{p}(n^{-1/2}\lambda^{-1/2r})$
and $o_{p}(\{\lambda{J}(\hat{\eta}^*-\eta^{*})
(V+\lambda{J})(\hat{\eta}^*-\eta_{0})\}^{1/2})$ by similar arguments for
(\ref{gn1_hzd:8}), and the terms in the third bracket is of
the order
\[
\{(V+\lambda{J})(\hat{\eta}-\eta^{*})\}^{1/2}
o_{p}\bigl(\{\lambda{J}(\hat{\eta}^*-\eta^{*})\}^{1/2}\bigr)
\]
by condition 3, Lemma \ref{vv_hzd:t} and (\ref{key_hzd:t}).
Putting all these together,
one obtains $(V+\lambda{J})(\hat{\eta}-\eta^{*})
=O_{p}(n^{-1}\lambda^{-1/r}+\lambda)$ and hence
$(V+\lambda{J})(\hat{\eta}-\eta_0)
=O_{p}(n^{-1}\lambda^{-1/r}+\lambda)$.
And an application of condition \hyperlink{con:lam}{A7} yields the final convergence
rates.\qed
\end{Step}
\noqed\end{pf*}
\begin{pf*}{Proof for the asymptotic properties of $\hat{\bolds{\beta}}$}
Let $P_n$ be the empirical measure of $(X_i,\Delta_i=1,Z_i),i=1,\ldots
,n$ such that
it is related to the empirical measure $Q_{n}$ of $(X_i,\Delta
_i,Z_i),i=1,\ldots,n$ by
$P_n f= \int f\,dP_n=\int\Delta f\,dQ_n =
n^{-1}\sum_{i=1}^{n}\Delta_i f(T_i,\Delta_i,Z_i)$.
Let $P$ be its corresponding\break (sub)probability measure.
Let $L_2(P)=\{f\dvtx\int f^2\,dP<\infty\}$ and $\|\cdot\|_2$ be the usual
$L_2$-norm. For any subclass $\mathcal{F}$ of $L_2(P)$ and any
$\varepsilon>0$, let $\mathcal{N}_{[\cdot]}(\varepsilon,\mathcal{F},L_2(P))$
be the bracketing number and
$J_{[\cdot]}(\delta,\mathcal{F},L_2(P))=\int_0^\delta
\sqrt{1+\log\mathcal{N}_{[\cdot]}(\varepsilon,\mathcal{F},L_2(P))}\,
d\varepsilon$.
\begin{lemma}\label{lem:entr}
Let
$m_0(t,u,w;\bolds{\beta},\eta)=u^{T}\bolds{\beta}+\eta(w)-\log s[\bolds
{\beta},\eta](t)$,
$m_1(t,u,\break w$; $s,\bolds{\beta},\eta)=1_{[s\le t]}\exp(u^{T}\bolds{\beta
}+\eta(w))$,
and
$m_2(t,u,w;s,\bolds{\beta},\eta,f)=1_{[s\le t]}f(u,\break w)\exp(u^{T}\bolds
{\beta}+\eta(w))$.
Define the classes of functions
\begin{eqnarray*}
\mathcal{M}_0(\delta)
&=&\{m_0\dvtx\|\bolds{\beta}-\bolds{\beta}_0\|\le\delta,\|\eta-\eta_0\|_2\le
\delta\},\\
\mathcal{M}_1(\delta)
&=&\{m_1\dvtx s\in\mathcal{T},\|\bolds{\beta}-\bolds{\beta}_0\|\le\delta,\|
\eta-\eta_0\|_2\le\delta\},\\
\mathcal{M}_2(\delta)
&=&\{m_2\dvtx s\in\mathcal{T},\|\bolds{\beta}-\bolds{\beta}_0\|\le\delta,\|
\eta-\eta_0\|_2\le\delta,
\|h\|_2\le\delta\}.
\end{eqnarray*}
Then $J_{[\cdot]}(\delta,\mathcal{M}_0,L_2(P))\le c_0q_n^{1/2}\delta$ and
$J_{[\cdot]}(\delta,\mathcal{M}_j,L_2(P))=c_j\delta\{q_n+\log(1/\break\delta)\}^{1/2}$
$j=1,2$.
\end{lemma}
\begin{pf}
The proof is similar to that of Corollary A.1 in \cite{huang99} and
thus omitted
here.
\end{pf}
\begin{lemma}\label{lem:i2}
\[
\sup_{t\in\mathcal{T}} |
s^{-1}[\bolds{\beta}_0,\hat{\eta}](t)
s[U;\bolds{\beta}_0,\hat{\eta}](t)-s_n^{-1}[\bolds{\beta}_0,\hat{\eta}](t)
s_n[U;\bolds{\beta}_0,\hat{\eta}](t) |=o_p(n^{-1/2}).
\]
\end{lemma}
\begin{pf}
Write
\begin{eqnarray*}
&&s^{-1}[\bolds{\beta}_0,\hat{\eta}](t)
s[U;\bolds{\beta}_0,\hat{\eta}](t)-s_n^{-1}[\bolds{\beta}_0,\hat{\eta}](t)
s_n[U;\bolds{\beta}_0,\hat{\eta}](t) \\
&&\qquad=\frac{s[\bolds{\beta}_0,\hat{\eta}](t)A_{1n}(t)-
s[U;\bolds{\beta}_0,\hat{\eta}](t)A_{2n}(t)}
{s[\bolds{\beta}_0,\hat{\eta}](t)s_n[\bolds{\beta}_0,\hat{\eta}](t)},
\end{eqnarray*}
where
$A_{1n}=s_n[U;\bolds{\beta}_0,\hat{\eta}](t)-s[U;\bolds{\beta}_0,\hat
{\eta}](t)$
and $A_{2n}=s_n[\bolds{\beta}_0,\hat{\eta}](t)-s[\bolds{\beta}_0,\break\hat
{\eta}](t)$.
Note that $q_n=o(n^{1/2})$, hence the result follows from
Lemma 3.4.2 of \cite{vaart96} and Lemma \ref{lem:entr}.
\end{pf}
\begin{pf*}{Proof of Theorem \ref{th:beta}}
Let $\gamma_n=n^{-1/2}$. To
prove \ref{th:beta}(i), we need to show that
$\forall\delta>0$, there exists a large constant $C$ such that
\[
P \Bigl\{\sup_{\|v\|=C}
\mathcal{L}_P(\bolds{\beta}_0+\gamma_n v)<\mathcal{L}_P(\bolds{\beta
}_0) \Bigr\}
\ge1-\delta.
\]
Consider $\mathcal{L}_P(\bolds{\beta}_0+\gamma_n v)-\mathcal{L}_P(\bolds
{\beta}_0)$. We can
decompose it to the sum of
$D_{n1}=l_p(\bolds{\beta}_0+\gamma_n v)-l_p(\bolds{\beta}_0)$ and the penalty
difference $D_{n2}$. As shown in \cite{scad01}, under the assumption of
$a_n=O(n^{-1/2})$ and $b_n=o(1)$, $n^{-1}D_{n2}$ is bounded by
%
\begin{equation}\label{eq:dn2}
\sqrt{s}\gamma_n a_n\|v\|+\gamma_n^2b_n\|v\|^2=C\gamma_n^2\bigl(\sqrt{s}+b_nC\bigr),
\end{equation}
where $s$ is the number of nonzero elements in $\bolds{\beta}_0$.

Applying the second order Taylor expansion to $n^{-1}D_{n1}$, gives
%
\begin{equation}\label{eq:dn1}
n^{-1}D_{n1}=\gamma_n\mathbf{v}^{T}\mathbf{J}_{1n}-\tfrac{1}{2}\gamma_n^2
\mathbf{v}^{T}\mathbf{J}_{2n}\mathbf{v}+o_p(n^{-1})
\end{equation}
with
$ 
\mathbf{J}_{1n}=P_n \{U-s_n^{-1}[\bolds{\beta}_0,\hat{\eta}](t)
s_n[U;\bolds{\beta}_0,\hat{\eta}](t) \},
$ 
and $\mathbf{J}_{2n}=P_n \{ s_n^{-2}[\bolds{\beta}_0,\hat{\eta}](t)\times\break [
s_n[UU^T;\bolds{\beta}_0,\hat{\eta}](t)
s_n[\bolds{\beta}_0,\hat{\eta}](t)-s_n[U;\bolds{\beta}_0,\hat{\eta}](t)
s_n[U;\bolds{\beta}_0,\hat{\eta}](t)^T
] \}$, where $U(u,\break w)\equiv u$.

Let $\bar{U}_n=\sum_{i=1}^{n}U_i/n$. Note that
$s_n^{-1}[\bolds{\beta}_0,\hat{\eta}](t)
s_n[\bar{U}_n;\bolds{\beta}_0,\hat{\eta}](t)=\bar{U}_n$. We have
%
\begin{eqnarray}\label{eq:j1n}
\mathbf{J}_{1n}&=&P_n \{U-\bar{U}_n-s_n^{-1}[\bolds{\beta}_0,\hat{\eta}](t)
s_n[U-\bar{U}_n;\bolds{\beta}_0,\hat{\eta}](t) \}\nonumber\\[-8pt]\\[-8pt]
&\equiv&I_{1n}+I_{2n}+I_{3n},\nonumber
\end{eqnarray}
where
\begin{eqnarray*}
I_{1n}&=& (P_n-P) \{U-\bar{U}_n-s^{-1}[\bolds{\beta}_0,\hat{\eta}](t)
s[U-\bar{U}_n;\bolds{\beta}_0,\hat{\eta}](t) \},\\[-1pt]
I_{2n}&=& P_n \{s^{-1}[\bolds{\beta}_0,\hat{\eta}](t)
s[U;\bolds{\beta}_0,\hat{\eta}](t)-s_n^{-1}[\bolds{\beta}_0,\hat{\eta}](t)
s_n[U;\bolds{\beta}_0,\hat{\eta}](t) \},\\[-1pt]
I_{3n}&=& P \{U-\bar{U}_n-s^{-1}[\bolds{\beta}_0,\hat{\eta}](t)
s[U-\bar{U}_n;\bolds{\beta}_0,\hat{\eta}](t) \}.
\end{eqnarray*}
Lemma 3.4.2 of \cite{vaart96} and Lemma \ref{lem:entr} indicate
that $(P_n-P) \{s^{-1}[\bolds{\beta}_0,\hat{\eta}](t)
s[U-\bar{U}_n;\bolds{\beta}_0,\hat{\eta}](t) \}=o_p(n^{-1/2})$,
where the fact that $q_n=o_p(n^{-1/2})$ is used again.
Also $(P_n-P)\{U-\bar{U}_n\}=O_p(n^{-1/2})$ by the LLN. Hence,
we have $I_{1n}=O_p(n^{-1/2})$. Lemma \ref{lem:i2} gives
$I_{2n}=o_p(n^{-1/2})$. Finally, by the boundedness assumption
$I_{3n}=P \{s^{-1}[\bolds{\beta}_0,\hat{\eta}](t)
s[\bar{U}_n-U;\bolds{\beta}_0,\hat{\eta}](t) \}
=O_p(E|\bar{U}_n-U|)=\break O_p(n^{-1/2})$. Hence,
$\mathbf{J}_{1n}=O_p(n^{-1/2})$.
Also, $\mathbf{J}_{2n}$ converges to $V(U)>0$. Thus, when $C$ is
sufficiently large, the second term in (\ref{eq:dn1}) dominates
both terms in (\ref{eq:dn2}). Theorem \ref{th:beta}(i)
follows.

Next,\vspace*{1pt} we shall show the sparsity of $\hat{\bolds{\beta}}$.
It suffices to show that for any given $\bolds{\beta}_1$
satisfying\vspace*{1pt}
$\|\bolds{\beta}_1-\bolds{\beta}_{10}\|=O_p(n^{-1/2})$ and any
$j=s+1,\ldots,d$,
$\partial\mathcal{L}_P(\bolds{\beta})/\partial\beta_j>0$ for
$0<\beta_j<Cn^{-1/2}$ and
$\partial\mathcal{L}_P(\bolds{\beta})/\partial\beta_j<0$ for
$-Cn^{-1/2}<\beta_j<0$. For $\beta_j\ne0$ and $j=s+1,\ldots,d$,
\[
n^{-1}\partial\mathcal{L}_P(\bolds{\beta})/\partial\beta_j
=P_n \{U_j-s_n^{-1}[\bolds{\beta},\hat{\eta}](t)
s_n[U_j;\bolds{\beta},\hat{\eta}](t) \}
-p'_{\theta_j}(|\beta_j|)\operatorname{sgn}(\beta_j).
\]
Similar to bounding $\mathbf{J}_{1n}$, the first term can be shown to be
$O_p(n^{-1/2})$. Recall that $\theta_j^{-1}=o(n^{1/2})$ and
$\lim\inf_{n\rightarrow\infty}\lim\inf_{u\rightarrow0^{+}}
\theta_{j}^{-1}p'_{\theta_{j}}(u)>0$. Hence, the sign of
$\partial\mathcal{L}_P(\bolds{\beta})/\partial\beta_j$ is completely
determined by that of $\beta_j$. Then $\hat{\bolds{\beta}}_2=\mathbf{0}$.

Lastly, we show the asymptotic normality of $\hat{\bolds{\beta}}_1$ using
the result in \cite{newey94}. Let $z_i=(X_i,\Delta_i,Z_i)$.
Note that $\hat{\bolds{\beta}}_1$ is the
solution of the estimating equation
%
\begin{equation}\label{eq:ee}
\sum_{i=1}^{n}M(z_i,\bolds{\beta}_1,\hat{\eta})-n\bolds{\zeta}_1=\mathbf{0},
\end{equation}
where $M(z,\bolds{\beta}_1,\eta)=\int\{U_1-s_n^{-1}[\bolds{\beta
}_1,\eta](t)
s_n[U_1;\bolds{\beta}_1,\eta](t) \}\,dN(t)$ and
$\bolds{\zeta}_1=\break (p'_{\theta_1}(|\beta_1|)\operatorname{sgn}(\beta_1),
\ldots,p'_{\theta_s}(|\beta_s|)\operatorname{sgn}(\beta_s) )^T$.
Let
\[
D(z,h)=\int\biggl\{\frac{s_n[U_1h;\bolds{\beta}_{10},\eta_0](t)}
{s_n[U_1h;\bolds{\beta}_{10},\eta_0](t)}-
\frac{s_n[U_1;\bolds{\beta}_{10},\eta_0](t)}{s_n[\bolds{\beta}_{10},\eta_0](t)}
\frac{s_n[h;\bolds{\beta}_{10},\eta_0](t)}{s_n[\bolds{\beta}_{10},\eta_0](t)}
\biggr\}\,dN(t)
\]
be the Fr\'echet derivative of $M(z,\bolds{\beta}_{10},\eta)$ at $\eta
_0$. Since the
convergence rate of $\hat{\eta}$ is $n^{-r/[2(r+1)]}=o(n^{-1/4})$,
the linearization assumption (Assumption 5.1) in \cite{newey94} is satisfied.
A derivation similar to bounding (\ref{eq:j1n}) can verify the
stochastic assumption
(Assumption 5.2) in \cite{newey94}. Direct calculation yields
$E[D(z,\eta-\eta_0)]=0$ for $\eta$ close to $\eta_0$. Then the
mean-square continuity
assumption (Assumption 5.3) in \cite{newey94} also holds with
$\alpha(z)\equiv0$. By Lemma 5.1 in \cite{newey94}, $\hat{\bolds{\beta
}}_1$ thus
has the same distribution as the solution to the equation
\[
\sum_{i=1}^{n}M(z_i,\bolds{\beta}_1,\eta_0)-n\bolds{\zeta}_1=\mathbf{0}.
\]
A straightforward simplification yields the result.
\end{pf*}
\noqed\end{pf*}
\end{appendix}

\section*{Acknowledgments}

We would like to thank the Associate Editor and two referees for
their insightful comments that have improved the article.

\printaddresses

\end{document}